\documentclass[11pt,reqno]{amsart}
\usepackage[T1]{fontenc}
\usepackage[utf8]{inputenc}
\usepackage{lmodern,microtype}
\usepackage[margin=1in]{geometry}
\usepackage{amsmath,amssymb,amsthm,mathtools}
\usepackage{enumitem}
\usepackage[hidelinks]{hyperref}
\hypersetup{pdftitle={Period bounds for Bergman kernel variation and the Suita defect},pdfauthor={Fusheng Deng, Gang Huang, Kefeng Liu}}
\setlist{leftmargin=*,itemsep=3pt}
\numberwithin{equation}{section}
\newtheorem{maintheorem}{Theorem}

\newtheorem{theorem}{Theorem}[section]
\newtheorem{proposition}[theorem]{Proposition}
\newtheorem{lemma}[theorem]{Lemma}
\newtheorem{corollary}[theorem]{Corollary}
\theoremstyle{definition}

\theoremstyle{remark}
\newtheorem{remark}[theorem]{Remark}

\title[Periods and the Suita defect]{Period bounds for Bergman kernel variation and the Suita defect}
\author[F. Deng]{Fusheng Deng}
\address{Fusheng Deng: School of Mathematical Sciences, University of Chinese Academy of Sciences, Beijing 100049, P. R. China}
\email{fshdeng@ucas.ac.cn}
\author[G. Huang]{Gang Huang}
\address{Gang Huang: School of Mathematical Sciences, University of Chinese Academy of Sciences, Beijing 100049, P. R. China}
\email{huanggang21@mails.ucas.ac.cn}
\author[K. Liu]{Kefeng Liu}
\address{Kefeng Liu: Mathematical Science Research Center, Chongqing University of Technology, Chongqing 400054, China}
\email{kefengliu@cqut.edu.cn}
\date{}
\makeatletter
\@ifundefined{subjclassname@2020}{%
  \@namedef{subjclassname@2020}{\textup{2020} Mathematics Subject Classification}%
}{}
\makeatother
\subjclass[2020]{Primary 30H20; Secondary 30C40, 32A25}
\keywords{Bergman kernel, exact differentials, Green function, Suita inequality, jet interpolation}
\begin{document}
\begin{abstract}
The periods of holomorphic differentials give an explicit lower bound
for the curvature of Bergman evaluation matrices. We prove this bound
along the Green sublevels of a smooth plane domain. For scalar evaluation at the Green pole, equality at one regular level holds precisely when that
sublevel is simply connected. For complete jets at the Green pole,
the matrices have continuous first derivatives, including at critical
levels, and their logarithmic determinant defects are integrals of
nonnegative curvature densities. A further estimate bounds the
accumulated period energy by the determinant defect, with a constant
independent of the jet order. On a multiply connected domain, all
these defects vanish up to the first critical Green value and become
strictly positive above it.
\end{abstract}
\maketitle

\section{Introduction}
The Bergman kernel of a plane domain depends on both its size and its
topology. The latter dependence can be measured by the periods of
holomorphic differentials. In this paper we show how these periods
enter the variation of the kernel along Green sublevels. The resulting
estimate gives a quantitative form of strict convexity and a bound
for the period energy when the Suita inequality is nearly an equality.

Let $D\Subset\mathbb C$ be connected with $C^\infty$ boundary, and
fix $q\in D$. Normalize the positive Green function by
\[
g_D(z,q)=\log\frac1{|z-q|}-\log c_D(q)+o(1),
\qquad z\longrightarrow q,
\]
and put
\[
\Phi=-g_D(\cdot,q),\qquad D_t=\{\Phi<t\},\qquad -\infty<t<0.
\]
The pole belongs to every sublevel. A value $t$ is called regular
if $\nabla\Phi$ has no zero on $\partial D_t$.

Bergman spaces use Euclidean area and inner products linear in the
first variable. The exact Bergman space is
\[
A_e^2(\Omega)=\{f\in A^2(\Omega):f(z)\,dz=dH
\text{ for some }H\in\mathcal O(\Omega)\}.
\]
It consists of the functions whose integrals around all closed curves
vanish. Write $\Pi_e$ for its orthogonal projection and $K_\Omega^e$
for its reproducing kernel. Then
\[
K_\Omega(q,q)-K_\Omega^e(q,q)
=\|(I-\Pi_e)K_\Omega(\cdot,q)\|_2^2.
\]
Thus the difference between the full and exact kernels measures the
period component of evaluation at $q$.

\subsection{The curvature estimate}
Let $L_1,\ldots,L_d$ be independent linear combinations of point
evaluations and holomorphic derivatives at fixed interior points.
On an interval where these points lie in $D_t$, let $k_{j,t}$ be
their Bergman representers and set
\[
M_{ij}(t)=\langle k_{j,t},k_{i,t}\rangle,
\qquad
M^e_{ij}(t)=\langle\Pi_e k_{j,t},\Pi_e k_{i,t}\rangle.
\]
For Hermitian matrices, $A\succeq B$ means that $A-B$ is positive
semidefinite.

\begin{maintheorem}\label{main:matrix}
At every regular Green level,
\begin{equation}\label{intro:matrix}
M''\succeq2M'(M+M^e)^{-1}M'.
\end{equation}
In particular, $(\log\det M)''\ge0$. For evaluation at the Green
pole $q$, write $K_t=K_{D_t}(q,q)$ and $K_t^e=K_{D_t}^e(q,q)$.
Then
\begin{equation}\label{intro:scalar}
(\log K_t)''\ge
\frac{K_t-K_t^e}{K_t+K_t^e}\bigl((\log K_t)'\bigr)^2.
\end{equation}
Equality in this scalar inequality at a regular level holds if and
only if $D_t$ is simply connected.
\end{maintheorem}

The matrix estimate has the exact remainder
\eqref{mx:matrix-remainder}. Its scalar form is
\eqref{sp:sharp-floor}--\eqref{sp:scalar-remainder}, and
Theorem~\ref{sp:green-rigidity} gives the equality criterion.
On a multiply connected regular sublevel the scalar inequality is
strict. Moreover, $K_t-K_t^e>0$ and $K_t'<0$, so
$(\log K_t)''>0$. Remark~\ref{sp:minimum-interpretation}
distinguishes the minimum energy under a fixed interior derivative constraint, with boundary values left free, from the energy of the harmonic extension prescribed by the Green flow.

\subsection{Determinant defects and stability}
For $\ell\ge1$, take the complete normalized jets
$L_jf=f^{(j)}(q)/j!$, $0\le j<\ell$, and denote their matrix by
$M^{(\ell)}$. Define
\begin{equation}\label{intro:jet}
F_\ell(t)=\log
\frac{\pi^\ell e^{\ell(\ell+1)t}\det M^{(\ell)}(t)}
{\ell!\,c_D(q)^{\ell(\ell+1)}}.
\end{equation}
Since $c_{D_t}(q)=e^{-t}c_D(q)$, the case $\ell=1$ is
\[
F_1(t)=\log\frac{\pi K_{D_t}(q,q)}{c_{D_t}(q)^2},
\]
the logarithmic defect in the Suita inequality.

\begin{maintheorem}\label{main:density}
The matrix $M^{(\ell)}$ is $C^1$ on $(-\infty,0)$, and its
derivative is locally absolutely continuous. At regular values set
$E_\ell=(\log\det M^{(\ell)})''$. Then
\begin{equation}\label{intro:density}
F_\ell(t)=\int_{-\infty}^t(t-s)E_\ell(s)\,ds,
\qquad -\infty<t\le0.
\end{equation}
At a regular value, $E_\ell(t)=0$ if and only if $D_t$ is simply
connected.
\end{maintheorem}

The values of $E_\ell$ at critical parameters do not affect the
integral. At $t=0$ we use the full-domain limits. Formula
\eqref{intro:density} expresses the determinant defect as the
accumulated curvature of the sublevels. In particular, $F_\ell\ge0$.

For these jets define
\[
P_t=M_t^{-1/2}(M_t-M_t^e)M_t^{-1/2},
\qquad
\rho_\ell(t)=\operatorname{tr}\bigl(P_t(2I-P_t)^{-1}\bigr).
\]
The eigenvalues of $P_t$ lie in $[0,1]$. We set $\rho_\ell=0$
at critical values, which does not affect its integrals. For $\ell=1$,
$\rho_1=(K_t-K_t^e)/(K_t+K_t^e)$.

\begin{maintheorem}\label{main:stability}
If $F_\ell(t_0)\le\varepsilon$, where $t_0\le0$ and
$\varepsilon\ge0$, then
\begin{equation}\label{intro:stability}
\int_{-\infty}^{t_0-\delta}\rho_\ell(t)\,dt
\le\frac{\varepsilon}{4\delta-2\varepsilon},
\qquad \delta>\varepsilon/2.
\end{equation}
The constant is independent of $\ell$.
\end{maintheorem}

Here $F_\ell$ is the full logarithmic determinant defect, without
division by $\ell$. If $e_1,\ldots,e_\ell$ is an orthonormal basis
of the span of the jet representers, then
\[
\operatorname{tr}P_t=\sum_{j=1}^{\ell}\|(I-\Pi_e)e_j\|_2^2,
\qquad \tfrac12\operatorname{tr}P_t\le\rho_\ell\le\operatorname{tr}P_t.
\]
Thus the estimate controls the sum of their period energies. Its
proof uses the reciprocal-trace identity of
Proposition~\ref{extp:reciprocal}, which in particular gives
\[
\frac d{dt}\operatorname{tr}\bigl(M(-M')^{-1}\bigr)\ge\rho_\ell.
\]
We also locate the threshold for positivity. If $D$ is multiply
connected and $t_*$ is the first critical value of $\Phi$, then
\[
F_\ell(t)=F_\ell'(t)=0\quad(t\le t_*),\qquad
F_\ell(t)>0,\quad F_\ell'(t)>0\quad(t_*<t<0).
\]
The threshold is the same for every complete jet.

\subsection{Earlier work and the proof}
The full and exact kernels and their period representations belong
to classical kernel theory \cite{BergmanSchiffer,Saitoh}.
Domain variation was studied by Komatsu \cite{Komatsu} and by
Maitani--Yamaguchi \cite{MaitaniYamaguchi}. Berndtsson proved
general positivity theorems \cite{Berndtsson06,Berndtsson09}.
Wang's Theorem~5.2, equation~(5.18), gives variation formulas for
derivative kernel evaluations; his Lemma~5.1, equation~(5.15),
expresses the flat-weight remainder by harmonic projections
\cite{Wang}. Thus jet variation and its harmonic contribution
are part of the existing theory. Corollary~\ref{mx:known-positivity}
recovers the corresponding convexity conclusions for the planar
Green family considered here.

Our further estimate uses the conjugacy of periods to count the
nonexact component twice in the interior energy. The resulting
weighted evaluation matrix is $(M+M^e)/2$, which gives the bound
in Theorem~\ref{main:matrix}. We determine equality for the Green
boundary data and then use the reciprocal trace to prove
Theorem~\ref{main:stability}. These explicit estimates and the
one-level equality criterion are the principal assertions of the
paper.

Regularity across singular domain variations also has antecedents.
Hamano--Maitani--Yamaguchi prove $C^1$ regularity of principal
functions and harmonic spans under the case~(c1) hypothesis of their
Lemma~5.1; their case~(c2) admits counterexamples
\cite{HamanoMaitaniYamaguchi}. The continuation proved here concerns
complete jet Bergman matrices along Green sublevels. It uses the
boundary estimate $|k_{r,t}|\le C_r|\nabla\Phi|$, including at
critical points of higher multiplicity.

Suita formulated the capacity inequality \cite{Suita}; B{\l}ocki
proved it for plane domains \cite{Blocki}, and Guan--Zhou obtained
the equality characterization on open Riemann surfaces
\cite{GuanZhou}. Corollaries~\ref{mx:jet-corollary}
and~\ref{mx:rigidity} recover the determinant inequality and its
endpoint equality case. The inequality also follows from the
higher-order estimates of B{\l}ocki--Zwonek \cite{BlockiZwonek}
by successive Schur complements; the deduction is given below.
Our formulas retain the curvature density and the period contribution.
Corollary~\ref{mx:sublevel} recovers the sublevel estimate of
\cite{BlockiSublevel,BlockiZwonek15}.

Section~\ref{sec:scalar} derives the classical variation formulas
from the Dirichlet problem. Section~\ref{sec:period} gives the
period decomposition, the weighted evaluation lemma, and the scalar
bound with its equality criterion. Section~\ref{mx:section} proves
Theorem~\ref{main:matrix}. The critical-level estimates and the
determinant formulas are in Section~\ref{sec:critical}. The last
two sections prove stability and locate the positivity threshold.

\section{Variation at regular levels}
\label{sec:scalar}
We record the classical planar variation formulas in our normalization;
see \cite{Komatsu,MaitaniYamaguchi,Wang}. The proofs from the Dirichlet
problem fix the constants used in the period estimates.
Throughout this section the domains are bounded, connected, and have
$C^\infty$ boundary. We use $dA=dx\,dy$,
$\Delta=4\partial_z\partial_{\bar z}$, and the inner product
$\langle f,h\rangle=\int f\overline h\,dA$, linear in its first entry.
The positive Green function and Dirichlet Green kernel satisfy
\begin{equation}\label{eq:normalization}
g_\Omega(z,w)=\log\frac1{|z-w|}+O(1),\qquad
G_\Omega=\frac{g_\Omega}{2\pi},\qquad -\Delta_zG_\Omega=\delta_w.
\end{equation}
Write
\begin{equation}\label{eq:regular-part}
G_\Omega(z,w)=\frac1{2\pi}\log\frac1{|z-w|}+S_\Omega(z,w).
\end{equation}
The Green--Bergman identity \cite{GreenBergman} is
\begin{equation}\label{eq:green-bergman}
K_\Omega(z,w)=-4\partial_z\partial_{\bar w}G_\Omega(z,w)
\quad(z\ne w),\qquad
K_\Omega(q,q)=-4\partial_z\partial_{\bar w}S_\Omega(q,q).
\end{equation}
The derivatives are taken before setting $z=w=q$. The identity holds
for the full Bergman kernel, also on multiply connected domains.
Let $n$ be the outward unit normal and put
$P(x,\xi)=-\partial_{n_\xi}G(x,\xi)$. Write $\mathcal P f$
for the harmonic extension of boundary data $f$ and
$\Lambda f=\partial_n\mathcal P f$ for its normal derivative.
For complex data, Green's identity gives
\begin{equation}\label{eq:dtn-energy}
\int_{\partial\Omega}\overline f\,\Lambda f\,d\sigma
=\int_\Omega|\nabla\mathcal P f|^2\,dA.
\end{equation}
Since $G(z,\cdot)=0$ on the boundary,
$\partial_{\bar w}G(z,\xi)=-n_c(\xi)P(z,\xi)/2$, where
$n_c=n_x+in_y$. Therefore
\begin{equation}\label{eq:boundary-kernel}
K_\Omega(z,\xi)=2n_c(\xi)\partial_zP(z,\xi).
\end{equation}

\subsection{Hadamard variations}
Let $\Omega_t$ depend $C^2$ on $t$ in the smooth spatial topology.
Choose normal parametrizations $X_t'=Vn$ and put
\[
\kappa=\operatorname{div}n,\qquad
A_n=X_t''\cdot n,\qquad B=A_n-\kappa V^2.
\]
All parameter derivatives keep the interior arguments fixed. When
we differentiate in a kernel variable, the domain and its boundary
velocity remain fixed. Thus differentiation in the evaluation
variable does not move the pole defining the Green sublevels.

\begin{proposition}[First two Green variations]\label{prop:hadamard}
At a fixed parameter, write $P_x(\xi)=P(x,\xi)$ and suppress the
domain subscript. Then
\begin{align}
G'(x,y)&=\int_{\partial\Omega}VP_xP_y\,d\sigma,\label{eq:Gprime}\\
G''(x,y)&=-2\int_{\partial\Omega}P_xV\Lambda(VP_y)\,d\sigma
+\int_{\partial\Omega}P_xBP_y\,d\sigma.\label{eq:Gsecond}
\end{align}
\end{proposition}
\begin{proof}
Fix a base parameter $t_0$ and write $\Omega=\Omega_{t_0}$.
Choose $U_0\Subset U\Subset\Omega$, domain diffeomorphisms $F_t$
with $F_{t_0}=\mathrm{id}$ and $F_t|_U=\mathrm{id}$, and
$\chi\in C_c^\infty(U)$ equal to $1$ near $\overline{U_0}$.
For $y\in U_0$, put $\Gamma_y(x)=(2\pi)^{-1}\log(1/|x-y|)$
and $R_t=G_{\Omega_t}(\cdot,y)-\chi\Gamma_y$. Then
\[
-\Delta R_t=f_y:=2\nabla\chi\cdot\nabla\Gamma_y
+(\Delta\chi)\Gamma_y,\qquad R_t|_{\partial\Omega_t}=0,
\]
where $(x,y)\mapsto f_y(x)$ is smooth on $\Omega\times U_0$,
with $x$-support in a fixed compact subset of
$U\setminus\overline{U_0}$. Pulling back gives $L_tu_t=f_y$, where
\[
u_t=R_t\circ F_t,\qquad
L_tu=\bigl[-\Delta(u\circ F_t^{-1})\bigr]\circ F_t.
\]
For $m\ge0$ and $0<\alpha<1$, set
\[
E_m=\{v\in C^{m+2,\alpha}(\overline\Omega):
v|_{\partial\Omega}=0\}.
\]
The operators $L_t:E_m\to C^{m,\alpha}(\overline\Omega)$ depend
$C^2$ on $t$ in operator norm. Schauder theory makes
$L_{t_0}=-\Delta$ an isomorphism \cite[Chapter 6]{GT}; a Neumann
series therefore gives a $C^2$ family of inverses near $t_0$.
Since $L_t$ is independent of $y$,
\[
\partial_y^\beta u_t=L_t^{-1}\partial_y^\beta f_y
\]
is $C^2$ in $t$ with values in $E_m$ for every multi-index $\beta$.
As $m$ is arbitrary, returning to the original coordinates gives
$C^2$ dependence of $S_{\Omega_t}$ in every fixed mixed interior
derivative norm, and smooth boundary traces of the first two
Eulerian derivatives. The singularity is independent of $t$, so
$G'$ and $G''$ are harmonic through $y$; their diagonal values are
interpreted as the corresponding derivatives of $S_{\Omega_t}$.

Now let $h_t=G_{\Omega_t}(\cdot,y)$ and $h=h_{t_0}$.
Differentiating $h_t(X_t)=0$ once gives
$h'=-V\partial_nh=VP_y$ on the boundary, hence
$h'=\mathcal P(VP_y)$. A second differentiation gives
\[
h''+2V\partial_nh'+V^2\partial_{nn}h+A_n\partial_nh=0
\quad\text{on }\partial\Omega.
\]
For the unit tangent $T$, the boundary condition and curvature
convention give $D^2h[T,T]=\kappa\partial_nh$. Harmonicity therefore
implies $\partial_{nn}h=-\kappa\partial_nh$, and consequently
\[
h''=-2V\Lambda(VP_y)+BP_y
\quad\text{on }\partial\Omega.
\]
Poisson representation for $h'$ and $h''$ proves
\eqref{eq:Gprime}--\eqref{eq:Gsecond}.
\end{proof}

\begin{corollary}[Scalar variation]\label{thm:scalar}
Fix $q$ in the common interior and put
$K=K_\Omega(q,q)$, $k=K_\Omega(\cdot,q)$,
$b(\xi)=V(\xi)\left.\partial_{\bar w}P(w,\xi)\right|_{w=q}$,
and $u=\mathcal P b$. Then
\begin{align}
K'&=-4u_z(q)=-\int_{\partial\Omega}V|K_\Omega(q,\xi)|^2\,d\sigma,
\label{eq:Kprime}\\
K''&=16(\|u_z\|_2^2+\|u_{\bar z}\|_2^2)
-\int_{\partial\Omega}B|K_\Omega(q,\xi)|^2\,d\sigma,\label{eq:Ksecond}\\
KK''-(K')^2
&=-K\int_{\partial\Omega}B|K_\Omega(q,\xi)|^2\,d\sigma
+16K\|u_{\bar z}\|_2^2 \notag\\
&\quad+16K\left\|u_z-\frac{u_z(q)}Kk\right\|_2^2.
\label{eq:scalar-two-square}
\end{align}
\end{corollary}
\begin{proof}
Proposition~\ref{prop:hadamard} permits the first two parameter
derivatives to commute with $\partial_z\partial_{\bar w}$, with
diagonal values taken through $S_{\Omega_t}$. Apply
$-4\partial_z\partial_{\bar w}$ to
\eqref{eq:Gprime}--\eqref{eq:Gsecond}, then set $z=w=q$;
$V$, $B$, and $\Lambda$ remain fixed during these kernel derivatives.
Since $u=\mathcal P b$, the first variation gives $K'=-4u_z(q)$,
and \eqref{eq:boundary-kernel} yields the boundary integral in
\eqref{eq:Kprime}.

For the second variation, \eqref{eq:dtn-energy} gives
\[
8\int_{\partial\Omega}\overline b\,\Lambda b\,d\sigma
=8\int_\Omega|\nabla u|^2\,dA
=16\bigl(\|u_z\|_2^2+\|u_{\bar z}\|_2^2\bigr).
\]
 The remaining term is
$-\int_{\partial\Omega}B|K_\Omega(q,\xi)|^2\,d\sigma$ by
\eqref{eq:boundary-kernel}, proving \eqref{eq:Ksecond}.

Finally, harmonicity and boundary regularity give $u_z\in A^2(\Omega)$.
Orthogonal projection onto the span of $k$, with $\|k\|_2^2=K$,
gives
\[
\|u_z\|_2^2=\frac{|u_z(q)|^2}{K}
+\left\|u_z-\frac{u_z(q)}Kk\right\|_2^2.
\]
Substituting into \eqref{eq:Ksecond} and using
$u_z(q)=-K'/4\in\mathbb R$ gives \eqref{eq:scalar-two-square}.
\end{proof}

\subsection{Green sublevels}
Let $\Phi=-g_D(\cdot,q)$, $D_t=\{\Phi<t\}$, $t<0$.
At a regular level, $V=1/\Phi_n=1/|\nabla\Phi|$.
The constant boundary value gives $\Phi_{TT}=\kappa\Phi_n$;
harmonicity gives $\Phi_{nn}=-\kappa\Phi_n$. Therefore
\[
A_n= X_t''\cdot n = V\partial_nV=-\Phi_{nn}/\Phi_n^3=\kappa V^2,\qquad B=0.
\]
The boundary term in \eqref{eq:scalar-two-square} therefore vanishes,
and $\log K_t$ is convex on every interval of regular levels.
The next section estimates its curvature in terms of periods.

\section{Periods and scalar curvature}
\label{sec:period}
The exact part of a holomorphic differential has a single-valued
primitive. Its orthogonal complement is finite dimensional and is
determined by the periods around the holes. We recall this
decomposition before applying it to the curvature.

\subsection{The period decomposition}

Let $\Omega$ be a bounded connected plane domain with $C^\infty$
boundary. Write its boundary components as
$\Gamma_0,\Gamma_1,\ldots,\Gamma_N$, with $\Gamma_0$ the outer
component. For $1\le j\le N$, choose a simple counterclockwise loop
$\gamma_j\subset\Omega$ enclosing the complementary component bounded
by $\Gamma_j$ and no other complementary component. Define
\[
 \ell_j(f)=\int_{\gamma_j}f(z)\,dz,
 \qquad \mathbf P(f)=(\ell_1(f),\ldots,\ell_N(f))^T.
\]
The exact Bergman space is
\[
 A^2_e(\Omega)=\ker\mathbf P.
\]
Thus $f\in A^2_e(\Omega)$ precisely when $f\,dz$ has a single-valued
holomorphic primitive. Interior estimates for holomorphic functions show
that each $\ell_j$ is bounded on $A^2(\Omega)$, so $A^2_e(\Omega)$ is
closed. Denote its reproducing kernel by $K^e_\Omega$. This is the
classical reduced Bergman kernel; see \cite{BergmanSchiffer,Saitoh}.

Let $\omega_j$ be the harmonic measure of $\Gamma_j$: its boundary
values are $1$ on $\Gamma_j$ and $0$ on every other component. Put
\[
 C_{jk}=\int_\Omega\nabla\omega_j\cdot\nabla\omega_k\,dA,
 \qquad r_j=2i\partial_z\omega_j,
 \qquad r=(r_1,\ldots,r_N)^T.
\]
Let $\Pi_e$ be the orthogonal projection onto $A^2_e(\Omega)$ and put
$\Pi_n=I-\Pi_e$.

\begin{lemma}[Period representation]\label{lem:period}
The matrix $C$ is real symmetric and positive definite, and
\[
 \ell_j(f)=\langle f,r_j\rangle,
 \qquad \langle r_k,r_j\rangle=C_{jk}.
\]
Consequently,
\begin{align}
 \Pi_nf&=r^TC^{-1}\mathbf P(f),
 &\|\Pi_nf\|^2&=\mathbf P(f)^*C^{-1}\mathbf P(f),
 \label{eq:period-energy}\\
 K_\Omega(z,q)-K^e_\Omega(z,q)
 &=r(z)^TC^{-1}\overline{r(q)}.
 \label{eq:kernel-difference}
\end{align}
If $u$ is a single-valued complex harmonic function of finite Dirichlet
energy and
$a=u_z$, $v=\overline{u_{\bar z}}$, then
\begin{equation}\label{eq:period-conjugacy}
 \mathbf P(a)+\overline{\mathbf P(v)}=0,
 \qquad \|\Pi_na\|=\|\Pi_nv\|.
\end{equation}
When $N=0$, the period formulas are understood with empty vectors and
$\Pi_n=0$.
\end{lemma}

\begin{proof}
Fix $j$ and choose $\chi_j\in C^\infty(\overline\Omega)$ that is equal
to $1$ in a collar of $\Gamma_j$ and to $0$ in collars of the other
components. Since $\omega_j-\chi_j\in W^{1,2}_0(\Omega)$, approximation
by compactly supported smooth functions and $\bar\partial f=0$ imply
\[
 \int_\Omega f\partial_{\bar z}(\omega_j-\chi_j)\,dA=0
 \qquad(f\in A^2(\Omega)).
\]
The support of $\bar\partial\chi_j$ is compactly contained in $\Omega$.
Apply Stokes' theorem on an inner parallel domain containing this support.
Its inner boundary components have clockwise induced orientation, and
$\chi_j$ is constant near them. Cauchy's theorem then gives
\[
 \langle f,r_j\rangle
 =-2i\int_\Omega f\partial_{\bar z}\omega_j\,dA
 =-2i\int_\Omega f\partial_{\bar z}\chi_j\,dA
 =\int_{\gamma_j}f\,dz.
\]
In particular, no boundary trace of $f$ is needed.

Writing $r_j=\partial_y\omega_j+i\partial_x\omega_j$, we obtain
\[
 \langle r_k,r_j\rangle
 =C_{jk}+i\int_\Omega
   (\partial_x\omega_k\partial_y\omega_j
    -\partial_y\omega_k\partial_x\omega_j)\,dA.
\]
The last integral is
$\int_{\partial\Omega}\omega_k\,d\omega_j$, which is zero because
the harmonic measures are constant on each boundary component. This
proves the Gram identity.

For $c\in\mathbb R^N$, set $u=\sum_{j=1}^N c_j\omega_j$.
Then $c^TCc=\int_\Omega|\nabla u|^2\,dA\ge0$, with equality
only when $u$ is constant. Since $u|_{\Gamma_0}=0$ and
$u|_{\Gamma_j}=c_j$, equality forces $c=0$, proving that
$C$ is positive definite.

The orthogonal complement of the common kernel of the $\ell_j$ is the
span of their linearly independent Riesz representatives $r_j$. Write
$\Pi_nf=r^T\beta$. Taking periods gives $C\beta=\mathbf P(f)$, and therefore
\[
 \|\Pi_nf\|^2
 =\beta^*C\beta=\mathbf P(f)^*C^{-1}\mathbf P(f).
\]
For $k_q=K_\Omega(\cdot,q)$, the reproducing property gives
$\mathbf P(k_q)=\overline{r(q)}$. Also $\Pi_e k_q=K^e_\Omega(\cdot,q)$,
because both vectors represent evaluation at $q$ on $A^2_e(\Omega)$.
This proves \eqref{eq:kernel-difference}.

The functions $a$ and $v$ are holomorphic and square integrable.
Integrating
$du=a\,dz+\overline v\,d\bar z$ around each $\gamma_j$ gives the first
identity in \eqref{eq:period-conjugacy}. Since $C^{-1}$ is real,
\[
 \mathbf P(v)^*C^{-1}\mathbf P(v)
 =\mathbf P(a)^TC^{-1}\overline{\mathbf P(a)}
 =\mathbf P(a)^*C^{-1}\mathbf P(a).
\]
The last equality follows because the quadratic form is real. The
norm identity follows from \eqref{eq:period-energy}.
\end{proof}

\begin{corollary}[Strict positivity of the period part]
\label{prop:strict-period}
If $\Omega$ is multiply connected, then
\[
 K_\Omega(q,q)-K^e_\Omega(q,q)>0
 \qquad(q\in\Omega).
\]
\end{corollary}

\begin{proof}
For $k_q=K_\Omega(\cdot,q)$, orthogonal projection gives
$K_\Omega(q,q)-K^e_\Omega(q,q)=\|\Pi_nk_q\|^2\ge0$.
Suppose equality holds. Then $k_q\in A^2_e(\Omega)$, so there is a
single-valued holomorphic $F$ with $F'=k_q$ and $F(q)=0$.
In particular, $F'(q)=K_\Omega(q,q)>0$.
The Green--Bergman identity and boundary regularity imply that
$k_q$, and hence $F$, are smooth up to $\partial\Omega$.

For $z\ne q$, set
$\mathcal B(z)=\overline{F(z)}+(2/\pi)\partial_qg_\Omega(z,q)$.
The conjugate Green--Bergman identity and the logarithmic singularity
give, respectively,
\[
 \partial_{\bar z}\mathcal B
 =\overline{k_q}+\frac2\pi\partial_{\bar z}\partial_qg_\Omega=0,
 \qquad
 \mathcal B(z)=\frac1{\pi(z-q)}+O(1)\quad(z\to q).
\]
Thus $\mathcal B$ is holomorphic off $q$ with a simple pole there.
Since $g_\Omega(z,q)=0$ for every boundary point $z$ and interior pole
$q$, we have $\partial_qg_\Omega=0$ and hence
$\mathcal B=\overline F$ on $\partial\Omega$.

Because $F(q)=0$, the product $F\mathcal B$ extends holomorphically
across $q$ and has real boundary values $|F|^2$.
Its imaginary part vanishes by the maximum principle, so holomorphicity
gives
\[
 F\mathcal B\equiv C:=\frac{F'(q)}\pi
                   =\frac{K_\Omega(q,q)}\pi>0.
\]
Consequently, $F$ has only its simple zero at $q$, and $|F|^2=C$
on every boundary component.

Set $\phi=F/\sqrt C$. Then $|\phi|=1$ on $\partial\Omega$;
since $\phi$ is nonconstant, the maximum modulus principle gives
$|\phi|<1$ in $\Omega$.
For any $\zeta\in\mathbb D$, the inequality
$|(\phi-\zeta)-\phi|=|\zeta|<|\phi|$ holds on every boundary component.
Rouch\'e's theorem therefore shows that $\phi-\zeta$, like $\phi$,
has exactly one zero in $\Omega$, counted with multiplicity.
Thus every $\zeta\in\mathbb D$ has a unique simple preimage, and
$\phi:\Omega\to\mathbb D$ is biholomorphic, contradicting multiple
connectivity.
\end{proof}

\begin{remark}
On the annulus $a<|z|<1$, write $L=\log(1/a)$. Then
\[
 \omega(z)=-\frac{\log|z|}{L},\qquad
 r(z)=-\frac{i}{Lz},\qquad C=\frac{2\pi}{L}.
\]
Thus \eqref{eq:kernel-difference} becomes
\[
 K_\Omega(z,q)-K^e_\Omega(z,q)
   =\frac1{2\pi Lz\bar q},
\]
the reproducing kernel of the one-dimensional nonexact space spanned
by $z^{-1}$.
\end{remark}

\subsection{Weighted functional evaluation}
On a fixed domain $\Omega$ as above, give $A^2(\Omega)$ the
equivalent Hilbert norm
\[
\langle f,h\rangle_*
=\langle\Pi_e f,\Pi_e h\rangle
 +2\langle\Pi_n f,\Pi_n h\rangle,
\qquad
\|f\|_*^2=\|\Pi_e f\|^2+2\|\Pi_n f\|^2.
\]
For a single-valued complex harmonic function $u$ of finite
Dirichlet energy, put $a=u_z$, $v=\overline{u_{\bar z}}$, and
\[
\mathcal D_\Omega(u)=\|u_z\|^2+\|u_{\bar z}\|^2
=\frac12\int_\Omega|\nabla u|^2\,dA.
\]
Lemma~\ref{lem:period} gives $\|\Pi_n v\|=\|\Pi_n a\|$, hence
\begin{equation}\label{sp:weighted-energy}
\mathcal D_\Omega(u)=\|a\|_*^2+\|\Pi_e v\|^2.
\end{equation}
The period component of $a$ therefore occurs twice in the energy.
The following interpolation identity will be used for both a single
functional and a finite family.

\begin{lemma}[Weighted functional evaluation]
\label{sp:harmonic-evaluation}\label{mx:minimization}
Let $L_1,\ldots,L_d$ be independent bounded complex-linear
functionals on $A^2(\Omega)$, with $L_if=\langle f,h_i\rangle$.
Set $Lf=(L_1f,\ldots,L_df)^T$ and
\begin{equation}\label{sp:weighted-gram}
\widetilde h_i=\Pi_eh_i+\tfrac12\Pi_nh_i,
\qquad
G_{ij}=L_i\widetilde h_j
=\tfrac12\bigl(\langle h_j,h_i\rangle
 +\langle\Pi_eh_j,\Pi_eh_i\rangle\bigr).
\end{equation}
The matrix $G$ is positive definite. For $y\in\mathbb C^d$, put
\begin{equation}\label{mx:minimizer}
a_y=\sum_{j=1}^d\widetilde h_j(G^{-1}y)_j.
\end{equation}
If $u$ is single-valued and harmonic with finite Dirichlet energy,
$a=u_z$, $v=\overline{u_{\bar z}}$, and $y=La$, then
\begin{equation}\label{mx:minimization-identity}
\mathcal D_\Omega(u)
=y^*G^{-1}y+\|a-a_y\|_*^2+\|\Pi_e v\|^2.
\end{equation}
The minimum of $\mathcal D_\Omega(u)$ under $L(u_z)=y$ is
$y^*G^{-1}y$, and it is attained. In particular, for a nonzero
bounded functional $L(f)=\langle f,h\rangle$,
\begin{equation}\label{sp:functional-bound}
2|L(u_z)|^2\le
\bigl(\|h\|^2+\|\Pi_eh\|^2\bigr)\mathcal D_\Omega(u),
\end{equation}
with optimal coefficient.
\end{lemma}
\begin{proof}
Define $T:\mathbb C^d\to A^2(\Omega)$ by
$Tc=(\widetilde h_1,\ldots,\widetilde h_d)c$.
Since $L_if=\langle f,\widetilde h_i\rangle_*$, we have
$LT=G$ and $\|Tc\|_*^2=c^*Gc$. Independence of the $L_i$ makes
$T$ injective, hence $G$ positive definite.
Now $a_y=TG^{-1}y$ satisfies $La_y=y$. For $La=y$, the difference
$r=a-a_y$ satisfies $Lr=0$ and
$\langle r,a_y\rangle_*=(G^{-1}y)^*Lr=0$. Thus
\[
\|a\|_*^2=\|a_y\|_*^2+\|a-a_y\|_*^2
=y^*G^{-1}y+\|a-a_y\|_*^2,
\]
which, together with \eqref{sp:weighted-energy}, proves
\eqref{mx:minimization-identity}.

By Lemma~\ref{lem:period}, choose $v_n\in\operatorname{ran}\Pi_n$
with $\mathbf P(v_n)=-\overline{\mathbf P(a_y)}$.
The form $a_y\,dz+\overline{v_n}\,d\bar z$ is closed and has zero
periods, hence equals $du$ for a single-valued harmonic function
of finite Dirichlet energy. Then $u_z=a_y$ and
$\Pi_e\overline{u_{\bar z}}=0$, so both remainder terms vanish.
For $d=1$, $G=(\|h\|^2+\|\Pi_eh\|^2)/2$ gives
\eqref{sp:functional-bound}; attainment for $y\ne0$ proves sharpness.
\end{proof}

For point evaluation at $q$, write
\[
k=K_\Omega(\cdot,q),\quad k_e=\Pi_ek,\quad k_n=\Pi_nk,
\qquad K=K_\Omega(q,q),\quad K^e=K_\Omega^e(q,q).
\]
Then
\begin{equation}\label{sp:kernel-norms}
\|k_e\|^2=k_e(q)=K^e,\qquad
\|k_n\|^2=k_n(q)=K-K^e.
\end{equation}
The constant function is exact, so $K^e>0$. The lemma gives
\begin{equation}\label{sp:harmonic-bound}
2|u_z(q)|^2\le(K+K^e)\mathcal D_\Omega(u).
\end{equation}
The boundary trace is free in this minimum-energy problem.
We next consider the trace prescribed by the Green motion.

\subsection{The scalar curvature bound}\label{sp:section}
Return to a regular Green sublevel $D_t$. Use the notation
$K=K_t$, $K^e=K_t^e$, $k=K_{D_t}(\cdot,q)$, and
$k_e=\Pi_ek$, $k_n=\Pi_nk$. Let $u$ be the harmonic function
in Corollary~\ref{thm:scalar}. Since $B=0$, its variation
formulas are
\begin{equation}\label{sp:green-variations}
K'=-4u_z(q),\qquad K''=16\mathcal D_{D_t}(u).
\end{equation}
Put
\begin{equation}\label{sp:scalar-minimizer}
a=u_z,\qquad v=\overline{u_{\bar z}},\qquad
a_*=-\frac{K'}{2(K+K^e)}\left(k_e+\tfrac12k_n\right).
\end{equation}

\begin{theorem}[Scalar curvature decomposition]
\label{sp:sharp-curvature}
At every regular Green level,
\begin{align}
(\log K)''
&=\frac{K-K^e}{K+K^e}\bigl((\log K)'\bigr)^2+R_*(t),
\label{sp:sharp-floor}\\
R_*(t)&=\frac{16}{K}
\bigl(\|a-a_*\|_*^2+\|\Pi_ev\|^2\bigr)\ge0.
\label{sp:scalar-remainder}
\end{align}
Consequently,
\begin{equation}\label{sp:curvature-bound}
(\log K_t)''\ge
\frac{K_t-K_t^e}{K_t+K_t^e}\bigl((\log K_t)'\bigr)^2.
\end{equation}
Equality holds precisely when $a=a_*$ and $\Pi_ev=0$.
If $D_t$ is multiply connected, then $(\log K_t)''>0$.
\end{theorem}
\begin{proof}
Apply Lemma~\ref{sp:harmonic-evaluation} to point evaluation with
$y=u_z(q)=-K'/4$. Its weighted Gram matrix is $G=(K+K^e)/2$
and its minimizer is $a_*$. Hence
\[
\mathcal D_{D_t}(u)
=\frac{(K')^2}{8(K+K^e)}+\|a-a_*\|_*^2+\|\Pi_ev\|^2.
\]
Multiply by $16/K$, use \eqref{sp:green-variations}, and
subtract $(K'/K)^2$ to prove the decomposition. Nonnegativity
gives the inequality and its equality conditions.

If $D_t$ is multiply connected, Corollary~\ref{prop:strict-period}
gives $K-K^e>0$. Moreover,
\[
-K'=\int_{\partial D_t}
\frac{|K_{D_t}(q,\xi)|^2}{|\nabla\Phi(\xi)|}\,d\sigma>0.
\]
The weight is positive, and the kernel section cannot vanish on
the entire boundary by the maximum principle. This proves strict
positivity of the curvature.
\end{proof}

\subsection{Pointwise equality for the Green flow}
Fix a regular level and write $\Omega=D_t$. It is connected by
Lemma~\ref{lem:levels} and has finitely many smooth boundary
components. Differentiation in the evaluation variable in
$V\partial_{\bar q}P_\Omega(q,\cdot)$ keeps both $\Omega$ and
$V$ fixed.

\begin{lemma}[Rigidity of the Green boundary data]
\label{sp:boundary-rigidity}
Let $u$ be the harmonic extension in Corollary~\ref{thm:scalar}
for the Green speed $V=1/|\nabla\Phi|$ at a regular level.
If $\Pi_e\overline{u_{\bar z}}=0$, then $u$ is holomorphic.
\end{lemma}
\begin{proof}
Put $g=g_\Omega(\cdot,q)=g_D(\cdot,q)+t$ and
$k=K_\Omega(\cdot,q)$, so $V=|\nabla g|^{-1}$.
On the boundary, \eqref{eq:boundary-kernel} gives
$\partial_{\bar q}P_\Omega(q,\xi)=n_c(\xi)k(\xi)/2$, while
$g_z=-|\nabla g|\overline{n_c}/2$. The prescribed boundary value is
\begin{equation}\label{sp:boundary-quotient}
u=-\frac{k}{4g_z}\quad\text{on }\partial\Omega.
\end{equation}
The Hopf lemma ensures $g_z\ne0$ there.

Let $\Gamma_0$ be the outer boundary and
$\omega=(\omega_1,\ldots,\omega_N)^T$ the harmonic measures of
the inner components. The hypothesis and Lemma~\ref{lem:period}
give $\overline{u_{\bar z}}=2i(\partial_z\omega)^T\alpha$ for
some $\alpha\in\mathbb C^N$.
Set $B_0=-2i\omega^T\overline\alpha$; then
$(B_0)_{\bar z}=u_{\bar z}$ and $B_0|_{\Gamma_0}=0$.
Thus $F=u-B_0\in\mathcal O(\Omega)\cap C^\infty(\overline\Omega)$
and $F=u$ on $\Gamma_0$.

Since $g_z=-1/(2(z-q))+O(1)$ near $q$, the function
\[
J(z)=(z-q)\bigl(4g_z(z)F(z)+k(z)\bigr)
\]
extends holomorphically across $q$ and is smooth up to the boundary.
Equation~\eqref{sp:boundary-quotient} gives $J|_{\Gamma_0}=0$.
Let $\widetilde J$ be the zero extension of $J$ across an arc of
$\Gamma_0$. Since $J=0$ on the arc, integration by parts gives
$\bar\partial\widetilde J=0$ distributionally. Thus $\widetilde J$
is holomorphic and vanishes on an exterior open set; the identity
theorem gives $J\equiv0$ in $\Omega$.
Hence $4g_zF+k=0$ on $\partial\Omega$, and
\eqref{sp:boundary-quotient} yields $F=u$ on every component.
The harmonic function $B_0=u-F$ has zero boundary values, so
$B_0\equiv0$ and $u=F$ is holomorphic.
\end{proof}

\begin{theorem}[Pointwise rigidity of the Green remainder]
\label{sp:green-rigidity}
At every regular value $t<0$,
\begin{equation}\label{sp:rigidity}
R_*(t)=0\quad\Longleftrightarrow\quad D_t\text{ is simply connected}.
\end{equation}
On every multiply connected regular sublevel,
\[
(\log K_t)''>
\frac{K_t-K_t^e}{K_t+K_t^e}\bigl((\log K_t)'\bigr)^2,
\qquad
K_t''>\frac{2(K_t')^2}{K_t+K_t^e}.
\]
\end{theorem}
\begin{proof}
If $\Omega=D_t$ is simply connected, choose a Riemann map
$\psi:\Omega\to\mathbb D$ with $\psi(q)=0$. Then
\[
g=-\log|\psi|,\qquad
g_z=-\frac{\psi'}{2\psi},\qquad
k=\frac{\overline{\psi'(q)}}{\pi}\psi'.
\]
Equation~\eqref{sp:boundary-quotient} identifies the prescribed
harmonic extension as
$u=\overline{\psi'(q)}\psi/(2\pi)$. Thus $v=0$, $a=k/2$,
and \eqref{sp:green-variations} gives $K'=-2K$. Since $K^e=K$
and $k_n=0$, \eqref{sp:scalar-minimizer} gives $a_*=k/2=a$.
Hence $R_*=0$.

Conversely, $R_*=0$ gives $\Pi_ev=0$ and $a=a_*$.
Lemma~\ref{sp:boundary-rigidity} makes $u$ a single-valued
holomorphic function. Its derivative $a$ is therefore exact.
Projecting \eqref{sp:scalar-minimizer} onto the period space gives
\[
0=\Pi_na=-\frac{K'}{4(K+K^e)}k_n.
\]
Since $-K'>0$, one has $k_n=0$ and $K=K^e$.
Corollary~\ref{prop:strict-period} forces $D_t$ to be simply
connected. The strict inequalities follow from
\eqref{sp:sharp-floor} and
$(\log K)''=K''/K-(K'/K)^2$.
\end{proof}

\begin{remark}\label{sp:minimum-interpretation}
For a fixed regular sublevel $\Omega=D_t$, let $\mathcal H_\alpha$
be the single-valued complex harmonic functions of finite Dirichlet
energy satisfying $h_z(q)=\alpha$. Lemma~
\ref{sp:harmonic-evaluation} gives
\[
\min_{h\in\mathcal H_\alpha}\mathcal D_\Omega(h)
=\frac{2|\alpha|^2}{K+K^e}.
\]
This minimum is attained on multiply connected domains as well.
For the Green extension $u$, however,
\[
\mathcal D_\Omega(u)
-\min_{h\in\mathcal H_{u_z(q)}}\mathcal D_\Omega(h)
=\frac K{16}R_*(t).
\]
Theorem~\ref{sp:green-rigidity} characterizes precisely when this
gap vanishes. Fixing the entire boundary trace gives the ordinary
Dirichlet principle,
\[
\mathcal D_\Omega(u+\varphi)
=\mathcal D_\Omega(u)+\mathcal D_\Omega(\varphi),
\qquad \varphi\in W^{1,2}_0(\Omega;\mathbb C),
\]
for either topology. The strict gap above is pointwise; it does
not give a uniform positive lower bound under degeneration.
\end{remark}

\section{The matrix curvature estimate}\label{mx:section}

The weighted evaluation lemma applies to a finite family of evaluation
and derivative functionals. We first record the corresponding variation
formula, then apply the lemma to obtain the period bound and its exact
remainder. The variation formula is the planar case of the classical
kernel-vector formulas; compare \cite[Theorem~5.2]{Wang}.

\subsection{The matrix variation identity}

Let $\Omega_t$ be a $C^2$ deformation of bounded $C^\infty$ planar domains.
Assume that the boundary motion is normal, with velocity $V$, normal
acceleration $A_n$, and $B=A_n-\kappa V^2$.  Fix linearly independent
complex-linear functionals $L_1,\ldots,L_d$, each a finite linear combination
of point evaluations and holomorphic derivatives at fixed points in the common
interior.
Write $L_i f=\langle f,k_i\rangle_t$ and set
\[
 M_{ij}=L_i k_j=\langle k_j,k_i\rangle_t,
 \qquad S=\operatorname{span}\{k_1,\ldots,k_d\}.
\]
Let $\Pi_S$ be the orthogonal projection onto $S$. Then
$c^*Mc=\|\sum_j c_jk_j\|_2^2$, so $M$ is positive definite.
For a collection $h_1,\ldots,h_d$ in any complex Hilbert space, use the
convention
\[
 \operatorname{Gram}(h)_{ij}=\langle h_j,h_i\rangle.
\]
The representers $k_j$ extend smoothly to the boundary, because their
evaluation points stay in the interior.

For the Poisson kernel $P_t(z,\xi)$ define
\begin{equation}\label{mx:boundary-data}
 d_j(\xi)=L_j\bigl[z\longmapsto\partial_zP_t(z,\xi)\bigr],
 \quad b_j(\xi)=V(\xi)\overline{d_j(\xi)},
 \quad u_j=\mathcal P_t b_j,
 \quad a_j=\partial_z u_j,
 \quad \eta_j=\partial_{\bar z}u_j.
\end{equation}
Here $\mathcal P_t$ is harmonic extension.  The functions $a_j$ are
holomorphic and the functions $\eta_j$ are antiholomorphic.

Put $r_j=(I-\Pi_S)a_j$ and
\[
(J_B)_{ij}=\int_{\partial\Omega_t}B\,\overline{k_i}\,k_j\,d\sigma.
\]

\begin{corollary}[Matrix variation]\label{mx:two-square}
For the preceding deformation,
\begin{align}
M'_{ij}&=-4L_i a_j
=-\int_{\partial\Omega_t}V\,\overline{k_i}\,k_j\,d\sigma,
\label{mx:first}\\
Q:=M''-M'M^{-1}M'
&=16\operatorname{Gram}(\eta)+16\operatorname{Gram}(r)-J_B.
\label{mx:Q}
\end{align}
If $B\le0$, then $Q\succeq0$ and
\begin{equation}\label{mx:logdet}
(\log\det M)''=\operatorname{tr}(M^{-1}Q)\ge0.
\end{equation}
If $V>0$ on the boundary, then $-M'\succ0$.
\end{corollary}

\begin{proof}
Reproduction gives
$\overline{k_j(z)}=\langle K_t(\cdot,z),k_j\rangle_t
=L_j[K_t(\cdot,z)]$. By kernel symmetry,
\[
k_j(z)=\overline{L_j}^{(w)}K_t(z,w),\qquad
M_{ij}=L_i^{(z)}\overline{L_j}^{(w)}K_t(z,w),
\]
where $\overline{L_j}g=\overline{L_j(\overline g)}$ for antiholomorphic
$g$. Proposition~\ref{prop:hadamard} permits these fixed interior
operations to commute with the first two parameter derivatives,
with diagonal values taken through the regular part.
By \eqref{eq:green-bergman}, apply
$-4L_i^{(z)}\overline{L_j}^{(w)}\partial_z\partial_{\bar w}$ to
\eqref{eq:Gprime}--\eqref{eq:Gsecond}. Using $u_j=\mathcal P_t b_j$
and $\overline{k_i}=2n_cd_i$ from \eqref{eq:boundary-kernel}, we obtain
\[
\begin{aligned}
M'_{ij}&=-4\int_{\partial\Omega_t}Vd_i\overline{d_j}\,d\sigma
=-4L_i a_j,\\
M''_{ij}&=8\int_{\partial\Omega_t}\overline{b_i}\Lambda b_j\,d\sigma
-(J_B)_{ij}.
\end{aligned}
\]
The first identity proves \eqref{mx:first}. Polarizing
\eqref{eq:dtn-energy} gives
$\int_{\partial\Omega_t}\overline{b_i}\Lambda b_j\,d\sigma
=2(\langle a_j,a_i\rangle_t+\langle\eta_j,\eta_i\rangle_t)$, hence
\begin{equation}\label{mx:second}
M''=16\operatorname{Gram}(a)+16\operatorname{Gram}(\eta)-J_B.
\end{equation}

Put $C=(L_i a_j)=-M'/4=C^*$ and write
$(\Pi_Sa_1,\ldots,\Pi_Sa_d)=(k_1,\ldots,k_d)X$.
Since $L_i r_j=0$, we have $MX=C$. Orthogonality gives
\[
\operatorname{Gram}(a)
=C^*M^{-1}C+\operatorname{Gram}(r)
=\tfrac1{16}M'M^{-1}M'+\operatorname{Gram}(r).
\]
Substitution proves \eqref{mx:Q}. If $B\le0$, then $Q\succeq0$,
and Jacobi's formula gives
$(\log\det M)''=\operatorname{tr}(M^{-1}Q)\ge0$.
Finally, if $V>0$ and $c^*(-M')c=0$, \eqref{mx:first} implies
$\sum_jc_jk_j=0$ on $\partial\Omega_t$. The maximum principle
and independence of the representers give $c=0$.
\end{proof}

\subsection{The period contribution}
Fix a regular Green level, so $B=0$ and $V>0$, and put
\[
M^e_{ij}=\langle\Pi_e k_j,\Pi_e k_i\rangle,\qquad
\Sigma=M+M^e.
\]
Then $0\preceq M^e\preceq M$ and $\Sigma\succ0$.
For the family $(L_i)$, the weighted Gram matrix in
Lemma~\ref{sp:harmonic-evaluation} is $G=\Sigma/2$.

\begin{theorem}[Matrix period bound]\label{mx:sharper-floor}
At every regular Green level,
\begin{equation}\label{mx:matrix-floor}
M''\succeq2M'(M+M^e)^{-1}M'.
\end{equation}
Consequently,
\begin{equation}\label{mx:trace-floor}
(\log\det M)''\ge
\operatorname{tr}\!\left(
M^{-1}M'\bigl[2(M+M^e)^{-1}-M^{-1}\bigr]M'\right)\ge0.
\end{equation}
For $c\in\mathbb C^d$ write $u_c=\sum_jc_ju_j$,
$a_c=(u_c)_z$, $v_c=\overline{(u_c)_{\bar z}}$, and set
\begin{equation}\label{mx:weighted-minimizer}
a_{*,c}=-\frac12\sum_j
\left(\Pi_ek_j+\tfrac12\Pi_nk_j\right)
[\Sigma^{-1}M'c]_j.
\end{equation}
The exact remainder is
\begin{equation}\label{mx:matrix-remainder}
c^*\bigl(M''-2M'\Sigma^{-1}M'\bigr)c
=16\bigl(\|a_c-a_{*,c}\|_*^2+\|\Pi_ev_c\|^2\bigr).
\end{equation}
\end{theorem}
\begin{proof}
Fix $c\in\mathbb C^d$. Since $B=0$, \eqref{mx:first} and
\eqref{mx:second} give $La_c=-M'c/4$ and
$\mathcal D_\Omega(u_c)=c^*M''c/16$.
Apply \eqref{mx:minimization-identity} with $y=-M'c/4$ and
$G=\Sigma/2$. Then $a_y=a_{*,c}$ and
\[
\tfrac1{16}c^*M''c
=\tfrac18c^*M'\Sigma^{-1}M'c
+\|a_c-a_{*,c}\|_*^2+\|\Pi_ev_c\|^2.
\]
Multiplying by $16$ proves \eqref{mx:matrix-remainder}; nonnegativity
of the squared norms for every $c$ gives \eqref{mx:matrix-floor}.

Set $R=M^{-1/2}M^eM^{-1/2}$. Then $0\preceq R\preceq I$, and
$(I-R)(I+R)^{-1}$ has eigenvalues $(1-\lambda)/(1+\lambda)\ge0$
for $\lambda\in\operatorname{spec}(R)\subset[0,1]$. Thus
\begin{equation}\label{mx:bracket}
2(M+M^e)^{-1}-M^{-1}
=M^{-1/2}(I-R)(I+R)^{-1}M^{-1/2}\succeq0.
\end{equation}
Denote this matrix by $W$. Subtracting $M'M^{-1}M'$ from
\eqref{mx:matrix-floor} gives $Q\succeq M'WM'\succeq0$.
Taking traces after congruence by $M^{-1/2}$ and using
\eqref{mx:logdet} proves \eqref{mx:trace-floor}.
\end{proof}

Equality in \eqref{mx:matrix-floor} holds precisely when
$a_c=a_{*,c}$ and $\Pi_ev_c=0$ for every $c$. For scalar
evaluation at the Green pole, \eqref{mx:matrix-remainder} reduces
to \eqref{sp:sharp-floor}--\eqref{sp:scalar-remainder}.
Theorem~\ref{sp:green-rigidity} gives the scalar equality
criterion. This proves Theorem~\ref{main:matrix}.

The following records the convexity consequences in this Green
family of the positivity results of Maitani--Yamaguchi,
Berndtsson, and Wang
\cite{MaitaniYamaguchi,Berndtsson06,Berndtsson09,Wang}.

\begin{corollary}[Convexity along Green sublevels]
\label{mx:known-positivity}
Let $M$ be the evaluation matrix in Theorem~\ref{main:matrix}.
On every interval of regular Green levels containing the evaluation
points,
\[
Q=M''-M'M^{-1}M'\succeq0,
\qquad (\log\det M)''\ge0.
\]
For every fixed nonzero $c\in\mathbb C^d$,
\[
\frac{d^2}{dt^2}\log(c^*M(t)c)\ge0.
\]
In particular, $t\mapsto\log K_{D_t}(z,z)$ is convex for a fixed
interior point $z$ on such an interval.
\end{corollary}
\begin{proof}
Since $M^e\preceq M$, inverse order gives
$2(M+M^e)^{-1}\succeq M^{-1}$. Theorem~\ref{main:matrix}
therefore gives $Q\succeq0$, and taking the trace against $M^{-1}$
gives determinant convexity. For $f=c^*Mc$, Cauchy--Schwarz yields
\[
f''\ge c^*M'M^{-1}M'c
\ge\frac{|c^*M'c|^2}{c^*Mc}
=\frac{(f')^2}{f}.
\]
Thus $(\log f)''\ge0$. Taking one point evaluation proves the
last assertion.
\end{proof}

\section{Critical levels and determinant defects}
\label{sec:critical}
The regular-level formulas do not by themselves control what happens
when a sublevel changes topology. We prove that the jet matrices are
$C^1$ and that their first derivatives are locally absolutely continuous. These
facts exclude singular curvature mass and permit integration from
the small sublevels near the pole.

\subsection{Local Green estimates}
We use Green-function monotonicity and exhaustion convergence
\cite{Ransford}, together with Bergman kernel convergence under increasing
exhaustions; see \cite[Theorem~4]{Ramadanov}. Both convergences are locally uniform
away from the Green pole. Interior estimates give convergence of all
fixed derivatives after the logarithmic singularity has been subtracted.

\begin{lemma}[Connectedness and compact critical strips]\label{lem:levels}
For bounded smooth $D$, each $D_t$ is connected, and only finitely
many critical points lie in a fixed strip $a\le\Phi\le b<0$.
\end{lemma}
\begin{proof}
A component not containing $q$ is relatively compact away from
the pole, with harmonic $\Phi$, boundary value $t$, and values
less than $t$ inside, contradicting the minimum principle.
The strip is compact and avoids the pole. The holomorphic function
$\Phi_z$ is not identically zero, so its zeros there are finite.
\end{proof}

\begin{lemma}[Uniform boundary estimate]\label{lem:boundary-control}
For regular $t$ in a compact interval in $(-\infty,0)$,
\begin{equation}\label{eq:boundary-control}
|K_{D_t}(q,\xi)|\le C|\nabla\Phi(\xi)|,\quad \xi\in\partial D_t.
\end{equation}
The same bound holds for each fixed holomorphic derivative in the
first kernel variable at $q$, with an order-dependent constant.
\end{lemma}
\begin{proof}
Let $[a,b]\Subset(-\infty,0)$ be the parameter interval and choose
$B(q,2r)\Subset D_a$. For regular $t\in[a,b]$ and
$\xi\in\partial D_t$, the function $P_t(\cdot,\xi)$ is positive
harmonic on this fixed disk. Harnack's inequality and interior
estimates give $|\partial_z^{j+1}P_t(q,\xi)|\le C_jP_t(q,\xi)$,
with $C_j$ independent of $t,\xi$.
The identity $g_{D_t}(\cdot,q)=g_D(\cdot,q)+t$ gives
$P_t(q,\xi)=|\nabla\Phi(\xi)|/(2\pi)$. Differentiating
\eqref{eq:boundary-kernel} in its first variable therefore yields
\[
|\partial_z^jK_{D_t}(q,\xi)|
=2|\partial_z^{j+1}P_t(q,\xi)|
\le C_j|\nabla\Phi(\xi)|,\qquad j\ge0.
\]
The case $j=0$ proves \eqref{eq:boundary-control}, and $j\ge1$
gives the asserted derivative bounds.
\end{proof}

\begin{lemma}[Green comparison at a critical value]\label{lem:critical-green}
Let $c$ be the only critical value in $[c-2\eta,c+2\eta]$.
Choose disks $q\in B_2\Subset B_1\Subset B_0\Subset D_{c-2\eta}$.
For every finite $m$, the normalized regular parts satisfy
\begin{equation}\label{eq:regular-convergence}
\|S_t-S_c\|_{C^m(\overline B_2\times\overline B_2)}
\le C_m|t-c|
\end{equation}
at nearby regular $t$. On compact regular boundary arcs of $D_c$,
the Green kernels and their mixed derivatives converge smoothly
under the local level flow, with their other variable near $q$.
\end{lemma}
\begin{proof}
We first prove \eqref{eq:regular-convergence}, then establish smooth
convergence on regular boundary patches. Write $g_s=g_{D_s}$ and
$S_s=S_{D_s}$, and restrict all levels to $(c-\eta,c+\eta)$.
For regular $t$, uniformly in $w\in B_1$,
$g_t(z,w)\le g_D(z,w)\le C_0$ on $\partial B_0$, where
$t-\Phi(z)\ge\eta$. Comparison on
$D_t\setminus\overline B_0$ gives
\begin{equation}\label{eq:barrier}
0\le g_t(z,w)\le C(t-\Phi(z)),\qquad
z\in D_t\setminus\overline B_0,\quad w\in B_1.
\end{equation}
For regular $\sigma<t$, cancellation of the logarithmic singularities
makes $g_t-g_\sigma$ harmonic on $D_\sigma$, including at $w$.
Its boundary values lie between $0$ and $C(t-\sigma)$ by
\eqref{eq:barrier}, so the maximum principle gives
\begin{equation}\label{eq:lipschitz-green}
0\le g_t(z,w)-g_\sigma(z,w)\le C(t-\sigma),\qquad z\in D_\sigma.
\end{equation}
For $t>c$, exhaustion $D_\sigma\uparrow D_c$ in
\eqref{eq:lipschitz-green} gives $0\le g_t-g_c\le C(t-c)$ on $D_c$.
For regular $t<c<\rho$, monotonicity gives
$0\le g_c-g_t\le g_\rho-g_t\le C(\rho-t)$ on $D_t$;
let $\rho\downarrow c$. Thus
\[
|g_t(z,w)-g_c(z,w)|\le C|t-c|,\qquad
(z,w)\in(D_t\cap D_c)\times B_1.
\]

After division by $2\pi$, this bounds $S_t-S_c$ on $B_0\times B_1$.
The difference is locally bounded and separately harmonic, hence
$(\Delta_z+\Delta_w)(S_t-S_c)=0$ distributionally. Weyl's lemma
and interior estimates on $B_2\times B_2$ prove
\eqref{eq:regular-convergence}.

Near a compact regular boundary arc choose $U'\Subset U$, with
$\overline U$ disjoint from $\overline B_0$ and the critical points.
Let $\Xi_t$ be the time-$(t-c)$ flow of
$\nabla\Phi/|\nabla\Phi|^2$, so
$\Phi\circ\Xi_t=\Phi+t-c$ and $\Xi_t\to\mathrm{id}$ smoothly.
Set $U_+=U\cap D_c$, $U'_+=U'\cap D_c$, and
$d_t(x,w)=g_t(\Xi_t(x),w)-g_c(x,w)$.
The barrier \eqref{eq:barrier}, also valid at $c$ by exhaustion,
gives $|d_t(x,w)|\le C(c-\Phi(x))$.
On $c-\Phi\ge\delta>0$, the common-domain comparison and
$\Xi_t\to\mathrm{id}$ give uniform convergence.
Letting $t\to c$ and then $\delta\downarrow0$ yields $d_t\to0$
uniformly on $\overline{U_+}\times B_1$.

Put $L_tu=[\Delta(u\circ\Xi_t^{-1})]\circ\Xi_t$.
Then $L_t\to\Delta$ smoothly and
\[
L_td_t=F_t:=-(L_t-\Delta)g_c\quad\text{in }U_+,
\qquad d_t|_{U\cap\partial D_c}=0.
\]
Local boundary regularity and interior estimates in $w$ make $g_c$
smooth on this patch. For each $m\ge0$ and $0<\alpha<1$,
$F_t\to0$ in the spatial $C^{m,\alpha}$ norm, uniformly for $w\in B_1$.
The local boundary Schauder estimate \cite[Chapter 6]{GT} gives
\[
\|d_t\|_{C^{m+2,\alpha}(\overline{U'_+})}
\le C_m\bigl(\|d_t\|_{C^0(\overline{U_+})}
+\|F_t\|_{C^{m,\alpha}(\overline{U_+})}\bigr)\longrightarrow0.
\]
Here the norms are in $x$, $C_m$ is uniform in $t,w$, and
localization to $U'$ avoids $\partial U\cap D_c$.
Finally, $d_t(x,\cdot)$ is harmonic on $B_1$, so interior estimates
give $D_w^\beta d_t\to0$ uniformly on
$\overline{U_+}\times\overline B_2$.
Since $L_t$ is independent of $w$,
$L_tD_w^\beta d_t=D_w^\beta F_t$, whose right-hand side also tends
to zero in the required norms. The same boundary estimate proves
the asserted convergence of all fixed mixed derivatives.
\end{proof}

\subsection{Continuation of the jet matrices}
Fix $\ell\ge1$ and take
$L_rf=f^{(r)}(q)/r!$, $0\le r<\ell$. The corresponding matrix is
\begin{equation}\label{mx:jet-matrix}
(M_t^{(\ell)})_{rs}
=\frac{1}{r!s!}
(\partial_z^r\partial_{\bar w}^sK_{D_t})(q,q).
\end{equation}
Both the full and exact jet matrices are positive definite.
Indeed, polynomials interpolate arbitrary finite jets and belong
to the exact Bergman space of every bounded domain.

\begin{lemma}[Continuation through critical levels]\label{mx:critical}
The matrix $M_t^{(\ell)}$ is $C^1$ on $(-\infty,0)$ and its
derivative is locally absolutely continuous. The function
$\log\det M_t^{(\ell)}$ is convex, and its second distributional
derivative is an absolutely continuous measure. At regular values,
\begin{equation}\label{mx:jet-first}
(M_t^{(\ell)})'_{rs}
=-\int_{\partial D_t}
\frac{\overline{k_{r,t}(\xi)}k_{s,t}(\xi)}{|\nabla\Phi(\xi)|}
\,d\sigma(\xi).
\end{equation}
\end{lemma}
\begin{proof}
Write $M(t)=M_t^{(\ell)}$ and let $\mathcal C$ be the locally finite
set of critical values (Lemma~\ref{lem:levels}). At regular values,
\eqref{mx:jet-first} follows from \eqref{mx:first} with
$V=|\nabla\Phi|^{-1}$. Fix $c\in\mathcal C$.
Lemma~\ref{lem:critical-green} and the Green--Bergman identity give
$k_{r,t}\to k_{r,c}$ smoothly near $q$ and, under the local level
flow, on compact regular arcs of $\partial D_c$.
Consequently, $M_{rs}(t)=L_rk_{s,t}\to L_rk_{s,c}=M_{rs}(c)$,
and the integrands in \eqref{mx:jet-first} converge smoothly
on those arcs.

At each critical point $p_j$ on this level, choose $W_j$ holomorphic
with $\operatorname{Re}W_j=\Phi$. Factoring $W_j-W_j(p_j)$ and
taking a holomorphic root gives a conformal coordinate $\zeta_j$
such that $W_j-W_j(p_j)=\zeta_j^{m_j}$, with $m_j\ge2$.
Lemma~\ref{lem:boundary-control} and kernel symmetry give
\begin{equation}\label{mx:jet-boundary-estimate}
|k_{r,t}(\xi)|\le C_r|\nabla\Phi(\xi)|.
\end{equation}
Thus, for $U_{j,\delta}=\{|\zeta_j|<\delta\}$,
\begin{align*}
\left|\int_{\partial D_t\cap U_{j,\delta}}
\frac{\overline{k_{r,t}}k_{s,t}}{|\nabla\Phi|}\,d\sigma\right|
&\le C_{r,s}\int_{\partial D_t\cap U_{j,\delta}}
|\nabla\Phi|\,d\sigma\\
&=C_{r,s}\int_{\partial D_t\cap U_{j,\delta}}|dW_j|
\le C'_{r,s,j}\delta^{m_j}.
\end{align*}
Indeed, $W_j$ maps the level arcs to a vertical chord of length at
most $2\delta^{m_j}$, with multiplicity at most $m_j$.
The coordinate circles meet $\partial D_c$ transversely, so the
remaining arcs and their integrals converge from both sides.
The total omitted contribution is $O(\sum_j\delta^{m_j})$,
uniformly in $t$, since there are finitely many critical points.
Letting $t\to c$ and then $\delta\downarrow0$ gives a common finite
limit $M'(t)\to N_c$. Continuity of $M$ and the mean-value theorem
give $M'(c)=N_c$. Hence $M\in C^1(-\infty,0)$.

For $b\ne0$, put $f_b=b^*Mb>0$ and $h_b=\log f_b\in C^1$.
Corollary~\ref{mx:known-positivity} gives $h_b''\ge0$ at regular
values. On any compact interval $[u,v]\subset(-\infty,0)$,
integrate over its finitely many regular subintervals. Since
$h_b''\ge0$, monotone convergence and continuity of $h_b'$ at
the critical endpoints give
\[
h_b'(v)-h_b'(u)
=\int_{(u,v)\setminus\mathcal C}h_b''(t)\,dt\ge0.
\]
This identity proves convexity of $h_b$ and
$h_b'\in\mathrm{AC}_{\mathrm{loc}}$. Hence
$f_b'=f_bh_b'\in\mathrm{AC}_{\mathrm{loc}}$.
Polarization using the coordinate vectors $e_r$ and their sums
$e_r+e_s$, $e_r+ie_s$ gives $M'\in\mathrm{AC}_{\mathrm{loc}}$
entrywise.

Finally, $\psi=\log\det M$ satisfies
$\psi'=\operatorname{tr}(M^{-1}M')\in\mathrm{AC}_{\mathrm{loc}}$.
By \eqref{mx:logdet},
$\psi''=\operatorname{tr}(M^{-1}Q)\ge0$ almost everywhere.
Thus $\psi$ is convex and
$D^2\psi=\operatorname{tr}(M^{-1}Q)\,dt$ is an absolutely
continuous measure.
\end{proof}

\subsection{The determinant defect}

\begin{lemma}[Conformal jet determinant]\label{mx:conformal-jets}
If $\Omega$ is simply connected, bounded, and $q\in\Omega$, then
\begin{equation}\label{mx:conformal-determinant}
 \det M_\Omega^{(\ell)}
   =\frac{\ell!}{\pi^\ell}c_\Omega(q)^{\ell(\ell+1)}.
\end{equation}
For the Green family of a bounded smooth $D$, sufficiently deep levels
are simply connected, and therefore
\begin{equation}\label{mx:deep-level}
 \det M_t^{(\ell)}
 =\frac{\ell!}{\pi^\ell}c_D(q)^{\ell(\ell+1)}
    e^{-\ell(\ell+1)t}
 \quad\text{for all sufficiently negative }t.
\end{equation}
\end{lemma}

\begin{proof}
For the disk $B(0,R)$ its Bergman kernel expansion gives the diagonal
normalized-jet matrix
\[
 (M_{B(0,R)}^{(\ell)})_{rs}
   =\delta_{rs}\frac{r+1}{\pi R^{2r+2}},
 \qquad
 \det M_{B(0,R)}^{(\ell)}
   =\frac{\ell!}{\pi^\ell R^{\ell(\ell+1)}}.
\]
Let $h:\Omega\to\mathbb D$ be a Riemann map with $h(q)=0$.
Change of variables gives $\|U_hF\|_{A^2(\Omega)}=\|F\|_{A^2(\mathbb D)}$
for $U_hF=h'(F\circ h)$; its inverse is $U_{h^{-1}}$, so $U_h$ is unitary.
Let $L_\Omega$ and $L_{\mathbb D}$ collect the first $\ell$ normalized
jets at $q$ and $0$, respectively. The adjoint identity
$\langle f,L_\Omega^*e_s\rangle=L_sf
=\langle f,k_{s,\Omega}\rangle$ gives $L_\Omega^*e_s=k_{s,\Omega}$,
hence $M_\Omega^{(\ell)}=L_\Omega L_\Omega^*$, and likewise on $\mathbb D$.
Since $h(q)=0$, the $r$th Taylor coefficient of $U_hF$ is
$h'(q)^{r+1}F^{(r)}(0)/r!$ plus terms involving lower derivatives of $F$.
Thus $L_\Omega U_h=T_hL_{\mathbb D}$, where $T_h$ is lower triangular
with $(T_h)_{rr}=h'(q)^{r+1}$ for $0\le r<\ell$. Therefore
\[
\begin{aligned}
M_\Omega^{(\ell)}
&=L_\Omega L_\Omega^*
=L_\Omega U_hU_h^*L_\Omega^*\\
&=T_hL_{\mathbb D}L_{\mathbb D}^*T_h^*
=T_hM_{\mathbb D}^{(\ell)}T_h^*.
\end{aligned}
\]
Taking determinants gives
$\det M_\Omega^{(\ell)}=|\det T_h|^2\det M_{\mathbb D}^{(\ell)}
=|h'(q)|^{\ell(\ell+1)}\ell!/\pi^\ell$.
Since $g_\Omega(z,q)=\log(1/|h(z)|)$, comparison of the constant
terms at $q$ gives $c_\Omega(q)=|h'(q)|$, proving
\eqref{mx:conformal-determinant}.

For the Green family, write $\Phi(z)=\log|z-q|+H_0(z)$ near $q$,
with $H_0$ harmonic, and choose holomorphic $H$ with
$\operatorname{Re}H=H_0$. The local coordinate $h(z)=(z-q)e^{H(z)}$
satisfies $h'(q)=e^{H(q)}\ne0$ and $\Phi=\log|h|$.
Choose a neighborhood $U$ of $q$ on which $h$ is univalent,
and a disk $B(q,r)\Subset U$.
The minimum principle bounds $\Phi$ below on $D\setminus\overline{B(q,r)}$,
so $D_t\subset U$ for sufficiently negative $t$.
Decreasing $t$ further ensures $B(0,e^t)\subset h(U)$, and hence
$h(D_t)=h(U)\cap B(0,e^t)=B(0,e^t)$.
Thus $D_t$ is simply connected. Finally,
$g_{D_t}(\cdot,q)=g_D(\cdot,q)+t$ gives
$c_{D_t}(q)=e^{-t}c_D(q)$; substitution in
\eqref{mx:conformal-determinant} proves \eqref{mx:deep-level}.
\end{proof}

For the Green sublevels of $D$, recall the normalized defect
\begin{equation}\label{mx:F}
 F_\ell(t)=\log\left(
  \frac{\pi^\ell e^{\ell(\ell+1)t}\det M_t^{(\ell)}}
       {\ell!\,c_D(q)^{\ell(\ell+1)}}\right).
\end{equation}
At regular values put
\[
 E_\ell(t)=\operatorname{tr}\bigl((M_t^{(\ell)})^{-1}Q_t\bigr),
 \qquad Q_t=(M_t^{(\ell)})''
       -(M_t^{(\ell)})'(M_t^{(\ell)})^{-1}(M_t^{(\ell)})'.
\]
Give $E_\ell$ the value zero at critical values.

\begin{theorem}[Jet defect representation]\label{mx:jet-defect}
The function $F_\ell$ is $C^1$ on $(-\infty,0]$ (with the left
derivative at $0$), convex, nonnegative, and nondecreasing.
It vanishes on all sufficiently deep levels, and
\begin{align}
F_\ell'(t)&=\int_{-\infty}^tE_\ell(s)\,ds,
\label{mx:first-integral}\\
F_\ell(t)&=\int_{-\infty}^t(t-s)E_\ell(s)\,ds,
\qquad -\infty<t\le0.
\label{mx:second-integral}
\end{align}
In particular,
\begin{equation}\label{mx:endpoint-integral}
\log\frac{\pi^\ell\det M_D^{(\ell)}}
{\ell!\,c_D(q)^{\ell(\ell+1)}}
=\int_{-\infty}^0(-s)E_\ell(s)\,ds.
\end{equation}

On every regular level the density in these formulas is bounded below
by the right-hand side of \eqref{mx:trace-floor}, with the corresponding
jet matrices.
\end{theorem}

\begin{proof}
Put $M(t)=M_t^{(\ell)}$. Since
$F_\ell=\log\det M+\ell(\ell+1)t+C$ for a constant $C$,
Lemma~\ref{mx:critical} gives $F_\ell\in C^1$ and
$F_\ell'\in\mathrm{AC}_{\mathrm{loc}}$ on $(-\infty,0)$.
Corollary~\ref{mx:two-square} gives $F_\ell''=E_\ell\ge0$ at regular
values. Absolute continuity therefore yields
$D^2F_\ell=E_\ell(t)\,dt$ in distributions, so $F_\ell$ is convex.

By Lemma~\ref{mx:conformal-jets}, choose $t_0<0$ such that
$F_\ell=F_\ell'=0$ for $t\le t_0$; hence $E_\ell=0$ there.
For $t_0\le t<0$, the fundamental theorem of calculus and
Tonelli's theorem give
\[
\begin{aligned}
F_\ell'(t)&=\int_{t_0}^tE_\ell(s)\,ds,\\
F_\ell(t)&=\int_{t_0}^t\int_{t_0}^uE_\ell(s)\,ds\,du
=\int_{t_0}^t(t-s)E_\ell(s)\,ds.
\end{aligned}
\]
Since $E_\ell=0$ below $t_0$, these prove
\eqref{mx:first-integral}--\eqref{mx:second-integral} for $t<0$,
as well as $F_\ell'\ge0$ and $F_\ell\ge0$.

As $t\uparrow0$, Bergman kernel convergence under $D_t\uparrow D$
and interior Cauchy estimates give $M(t)\to M_D^{(\ell)}\succ0$,
hence $F_\ell(t)\to F_\ell(0)<\infty$.
Writing $(a)_+=\max\{a,0\}$, monotone convergence gives
\[
F_\ell(0)=\lim_{t\uparrow0}
\int_{-\infty}^0(t-s)_+E_\ell(s)\,ds
=\int_{-\infty}^0(-s)E_\ell(s)\,ds,
\]
which is \eqref{mx:endpoint-integral}.

Finally, Hopf's lemma gives $|\nabla\Phi|>0$ on $\partial D$, so the
level flow is smooth up to $t=0$. The regular variation formulas
give a finite limit of $M'(t)$, and therefore of
$F_\ell'(t)=\ell(\ell+1)+\operatorname{tr}(M^{-1}M')$.
Thus $F_\ell$ is $C^1$ up to $0$, and monotone convergence extends
\eqref{mx:first-integral} to this endpoint. The density bound follows
from Theorem~\ref{mx:sharper-floor} applied to these jet matrices.
\end{proof}

\begin{corollary}[Jet determinant inequality]\label{mx:jet-corollary}
If $D\Subset\mathbb C$ is a smooth domain and $q\in D$, then
\begin{equation}\label{mx:jet-Suita}
 \det M_D^{(\ell)}
    \geq\frac{\ell!}{\pi^\ell}c_D(q)^{\ell(\ell+1)}.
\end{equation}
\end{corollary}

\begin{proof}
The right-hand side of \eqref{mx:endpoint-integral} is nonnegative.
\end{proof}

The inequality also follows from the higher-order estimates of
B{\l}ocki--Zwonek \cite{BlockiZwonek}; see also
\cite[equation~(4.1) and Theorem~4.5]{XuZhou}. To see this,
let $s_r$ be the Schur complement of the $r$th normalized jet
after the lower jets. Orthogonal projection gives
\[
s_r=\sup\left\{\frac{|f^{(r)}(q)|^2}{(r!)^2}:
\|f\|_2\le1,\quad f^{(j)}(q)=0\ (0\le j<r)\right\}
\ge\frac{r+1}{\pi}c_D(q)^{2r+2}.
\]
Multiplying  these inequalities and using $\det M_D^{(\ell)}=\prod_{r=0}^{\ell-1}s_r$ proves
\eqref{mx:jet-Suita}. Equality in the product forces equality for
$r=0$; the ordinary Suita equality theorem \cite{GuanZhou} then
gives simple connectivity for bounded smooth $D$. The next
corollary also derives this conclusion from the curvature density.

\subsection{Equality and exhaustion}

\begin{corollary}[Equality]\label{mx:rigidity}
For complete jets at the Green pole and every regular $t<0$,
\[
 (\log\det M_t^{(\ell)})''=0
 \quad\Longleftrightarrow\quad D_t\text{ is simply connected}.
\]
For bounded smooth connected $D$, equality in \eqref{mx:jet-Suita} holds
if and only if $D$ is simply connected.
\end{corollary}

\begin{proof}
Fix a regular $t<0$ and write $M=M_t^{(\ell)}$.
Orthogonal projection gives $M-M^e\succeq0$.
If $D_t$ is multiply connected, Corollary~\ref{prop:strict-period}
gives $(M-M^e)_{00}=K_{D_t}(q,q)-K_{D_t}^e(q,q)>0$, so $M-M^e\ne0$.
For $R=M^{-1/2}M^eM^{-1/2}$, we have $0\preceq R\preceq I$ and
$R\ne I$; hence at least one eigenvalue of $R$ is less than $1$.
Formula~\eqref{mx:bracket} therefore gives
$W:=2(M+M^e)^{-1}-M^{-1}\succeq0$ with $W\ne0$.
By Corollary~\ref{mx:two-square}, $-M'\succ0$, so
$X=M'M^{-1/2}$ is invertible. Thus $X^*WX$ is nonzero and positive
semidefinite, and \eqref{mx:trace-floor} yields
\[
E_\ell(t)\ge\operatorname{tr}(M^{-1}M'WM')
=\operatorname{tr}(X^*WX)>0.
\]

If $D_t$ is simply connected, the regular level flow identifies
$D_s$ with $D_t$ for $s$ near $t$, so these $D_s$ are simply connected.
Formula~\eqref{mx:conformal-determinant} and
$c_{D_s}(q)=e^{-s}c_D(q)$ then give
$\log\det M_s^{(\ell)}=C-\ell(\ell+1)s$, with $C$ independent of $s$.
Hence $E_\ell(t)=0$.

For simply connected $D$, \eqref{mx:conformal-determinant} gives
equality in \eqref{mx:jet-Suita}. If $D$ is multiply connected,
Hopf's lemma and compactness give $|\nabla\Phi|>0$ in a boundary
collar. Thus, for some $\varepsilon>0$, every $s\in(-\varepsilon,0)$
is regular and $D_s$ is diffeomorphic to $D$.
The first part gives $E_\ell(s)>0$ there, so
\eqref{mx:endpoint-integral} implies
\[
F_\ell(0)=\int_{-\infty}^0(-s)E_\ell(s)\,ds
\ge\int_{-\varepsilon}^0(-s)E_\ell(s)\,ds>0.
\]
This proves strict inequality in \eqref{mx:jet-Suita} for multiply
connected $D$.
\end{proof}

Lemma~\ref{mx:critical}, Theorem~\ref{mx:jet-defect}, and
Corollary~\ref{mx:rigidity} prove Theorem~\ref{main:density}.

\begin{corollary}[Arbitrary Greenian domains]\label{mx:arbitrary}
Let $D\subset\mathbb C$ be any Greenian domain and $q\in D$.
The determinant inequality \eqref{mx:jet-Suita} holds.  The function
$F_\ell$ in \eqref{mx:F}, defined using its Green sublevels, is finite,
nonnegative, nondecreasing, and convex.
\end{corollary}

\begin{proof}
Choose a smooth bounded connected exhaustion $E_j\uparrow D$ with
$q\in E_1$, and set $m=\ell(\ell+1)$.
Green exhaustion gives $g_{E_j}\uparrow g_D$; subtracting the common
logarithmic singularity near $q$ yields $c_{E_j}(q)\downarrow c_D(q)>0$.
For $t<0$, put $E_{j,t}=\{z\in E_j:-g_{E_j}(z,q)<t\}$, and set
$E_{j,0}=E_j$, $D_0=D$. For every fixed $t\le0$, we have
$E_{j,t}\uparrow D_t$. Bergman kernel convergence and interior Cauchy
estimates therefore give $M_{E_{j,t}}^{(\ell)}\to M_{D_t}^{(\ell)}$
entrywise.

By Theorem~\ref{mx:jet-defect}, each $F_{\ell,E_j}$ is nonnegative,
nondecreasing, and convex. In particular,
\[
\det M_{E_{j,t}}^{(\ell)}
\ge\frac{\ell!}{\pi^\ell}e^{-mt}c_{E_j}(q)^m
\ge\frac{\ell!}{\pi^\ell}e^{-mt}c_D(q)^m>0.
\]
The limiting determinant is thus finite and strictly positive.
Together with $c_{E_j}(q)\to c_D(q)$, this permits passage to the
logarithm in \eqref{mx:F}, giving
$F_{\ell,E_j}(t)\to F_{\ell,D}(t)$ with a finite limit.
Nonnegativity, monotonicity, and convexity pass to this pointwise
limit. At $t=0$, the limiting determinant bound is exactly
\eqref{mx:jet-Suita}.
\end{proof}

\begin{corollary}[Sublevel kernel estimate]\label{mx:sublevel}
For a Greenian plane domain $D$, $q\in D$, and $t<0$,
\[
K_D(q,q)\ge e^{2t}K_{D_t}(q,q)\ge\frac{e^{2t}}{|D_t|}.
\]
The last expression is zero when $|D_t|=\infty$.
\end{corollary}
\begin{proof}
The first inequality is monotonicity of $F_1$. If $|D_t|<\infty$,
testing the Bergman extremal problem with $|D_t|^{-1/2}$ gives
the second inequality. The infinite-area case is immediate.
\end{proof}
This recovers the estimate in \cite{BlockiSublevel,BlockiZwonek15}.
Exhaustion gives convexity and the determinant inequality on general
Greenian domains. Local absolute continuity of the curvature measure
is asserted here only for bounded smooth domains.

\section{Period stability for complete jets}
\label{extp:section}

We now prove Theorem~\ref{main:stability}. The matrix inequality
controls a reciprocal trace. Comparing that trace with the
derivative of the determinant defect gives a bound independent
of the jet order.

Let $D\Subset\mathbb C$ be connected with $C^\infty$ boundary,
let $q\in D$, and let $D_t=\{-g_D(\cdot,q)<t\}$, $t<0$.
Fix $\ell\ge1$ and use the complete normalized jets
$L_jf=f^{(j)}(q)/j!$, $0\le j<\ell$.
In this section write
\[
 M=M_t^{(\ell)},\qquad M^e=(M_t^{(\ell)})^e,\qquad
 H=-M',\qquad Q=M''-M'M^{-1}M'
\]
and let $F_\ell$ be the normalized determinant defect in
\eqref{mx:F}. This is the full logarithmic defect, without division
by $\ell$. Put $m_\ell=\ell(\ell+1)$.
Here $H$ is defined for all $t<0$; expressions involving $M''$
are understood at regular values, hence almost everywhere.
At regular values define the normalized period matrix and its trace
coefficient by
\begin{equation}\label{extp:period-matrix}
 P=M^{-1/2}(M-M^e)M^{-1/2},\qquad
 \rho_\ell=\operatorname{tr}\bigl(P(2I-P)^{-1}\bigr).
\end{equation}
If $\mu_1(t),\ldots,\mu_\ell(t)$ are the eigenvalues of $P_t$, then
\[
 0\le\mu_j(t)\le1,\qquad
 \rho_\ell(t)=\sum_{j=1}^\ell\frac{\mu_j(t)}{2-\mu_j(t)}.
\]
Define $\rho_\ell=0$ at critical values.  No continuity of $M^e$
across a critical level will be needed.

\subsection{The eigenvalue bounds}

\begin{lemma}[Eigenvalue bounds]\label{extp:speeds}
The matrix $H=-M'$ is positive definite. If
$\lambda_1(t)\le\cdots\le\lambda_\ell(t)$ are the eigenvalues of
$A=M^{-1/2}HM^{-1/2}$, then each $\lambda_j$ is nonincreasing and
\begin{equation}\label{extp:speed-bounds}
0<\lambda_j(t)\le2j,\qquad
F_\ell'(t)=\sum_{j=1}^\ell(2j-\lambda_j(t)).
\end{equation}
\end{lemma}
\begin{proof}
Fix $c\in\mathbb C^\ell\setminus\{0\}$ and set
$L_cf=\sum_{j=0}^{\ell-1}\overline{c_j}f^{(j)}(q)/j!$.
The reproducing property gives the extremal formula
\[
f_c(t):=c^*M(t)c
=\sup_{0\ne f\in A^2(D_t)}
  \frac{|L_cf|^2}{\displaystyle\int_{D_t}|f|^2\,dA}>0,
\]
and the supremum is attained. By Lemma~\ref{mx:critical},
$h_c=\log f_c$ is $C^1$. Corollary~\ref{mx:known-positivity} makes
$h_c$ convex on each regular interval; continuity of $h_c'$ extends
convexity across the locally finite critical values.

For $s<t$, choose a maximizer $u$ in the formula for $f_c(t)$.
Since $D_t\setminus D_s$ contains an open set, the identity theorem
gives $0<\|u\|_{L^2(D_s)}<\|u\|_{L^2(D_t)}$.
Restriction preserves $L_cu$, so
$f_c(s)\ge |L_cu|^2/\|u\|_{L^2(D_s)}^2>f_c(t)$.
For $t<\tau<0$, convexity gives
$h_c'(t)\le(h_c(\tau)-h_c(t))/(\tau-t)<0$.
Thus $c^*Hc=-f_c'=-f_ch_c'>0$, proving $H\succ0$.

The generalized Rayleigh quotient
$R_t(c)=c^*H(t)c/(c^*M(t)c)=-h_c'(t)$ is nonincreasing in $t$.
The substitution $x=M^{1/2}c$ identifies $R_t(c)$ with
$x^*Ax/\|x\|^2$, so the min--max principle gives
\begin{equation}\label{extp:parallel-frame}
 \lambda_j(t)=
 \min_{\substack{\mathcal V\subset\mathbb C^\ell\\\dim\mathcal V=j}}
 \max_{c\in\mathcal V\setminus\{0\}}R_t(c).
\end{equation}
The subspaces in this formula are independent of $t$; monotonicity
of each $R_t(c)$ therefore implies monotonicity of every $\lambda_j$.

For sufficiently negative $t$, the proof of
Lemma~\ref{mx:conformal-jets} provides a fixed coordinate
$h:D_t\to B(0,e^t)$. The disk matrix and jet transformation established
there give $M=T\Delta T^*$, where $T$ is invertible and independent
of $t$, and
$\Delta(t)=\operatorname{diag}((j/\pi)e^{-2jt})_{j=1}^\ell$.
Hence $H=T(-\Delta')T^*$ and $-\Delta_{jj}'=2j\Delta_{jj}$.
The identity $A-\lambda I=M^{-1/2}(H-\lambda M)M^{-1/2}$ identifies
the eigenvalues of $A$ with the generalized eigenvalues of $(H,M)$;
the common transformation $T$ cancels in
\[
\det(A-\lambda I)
=\frac{\det(H-\lambda M)}{\det M}
=\frac{\det(-\Delta'-\lambda\Delta)}{\det\Delta}
=\prod_{j=1}^\ell(2j-\lambda).
\]
Thus the ordered eigenvalues on these levels are exactly
$2,4,\ldots,2\ell$.
Monotonicity and $H\succ0$ now give $0<\lambda_j(t)\le2j$ for all
$t<0$. Finally,
$F_\ell'=m_\ell+\operatorname{tr}(M^{-1}M')
=m_\ell-\operatorname{tr}A=\sum_{j=1}^\ell(2j-\lambda_j)$.
\end{proof}

\subsection{The reciprocal trace}
Set
\[
U(t)=\operatorname{tr}(MH^{-1})=\sum_{j=1}^\ell\lambda_j(t)^{-1},
\qquad U_0=\frac12\sum_{j=1}^\ell\frac1j,
\qquad B_\ell(t)=\int_{-\infty}^t\rho_\ell(s)\,ds.
\]
At regular levels define the remainder
\[
\mathcal R=M''-2H(M+M^e)^{-1}H.
\]
It is positive semidefinite by Theorem~\ref{mx:sharper-floor}.

\begin{proposition}[Reciprocal-trace identity]\label{extp:reciprocal}
The function $U$ is locally absolutely continuous, and
\begin{equation}\label{extp:reciprocal-identity}
U'(t)=\rho_\ell(t)+
\operatorname{tr}(MH^{-1}\mathcal R H^{-1})
\quad\text{almost everywhere}.
\end{equation}
In particular,
\begin{equation}\label{extp:reciprocal-derivative}
U'\ge\rho_\ell\ge0,
\end{equation}
and $\rho_\ell$ is integrable on $(-\infty,t]$ for every $t\le0$,
with
\begin{equation}\label{extp:integrated-U}
B_\ell(t)\le U(t)-U_0.
\end{equation}
The endpoint uses one-sided limits.
\end{proposition}
\begin{proof}
Since $H\succ0$ and $M'$ is locally absolutely continuous, so are
$H^{-1}$ and $U$. Using $H'=-M''$ and the definition of $\mathcal R$,
we obtain almost everywhere
\begin{align*}
U'
&=-\ell+\operatorname{tr}(MH^{-1}M''H^{-1})\\
&=-\ell+2\operatorname{tr}\bigl(M(M+M^e)^{-1}\bigr)
 +\operatorname{tr}(MH^{-1}\mathcal R H^{-1}).
\end{align*}
Since $M+M^e=M^{1/2}(2I-P)M^{1/2}$, cyclicity gives
$-\ell+2\operatorname{tr}(M(M+M^e)^{-1})
=\operatorname{tr}(P(2I-P)^{-1})=\rho_\ell$, proving
\eqref{extp:reciprocal-identity}. Moreover,
$\operatorname{tr}(MH^{-1}\mathcal R H^{-1})
=\operatorname{tr}(\mathcal R^{1/2}H^{-1}MH^{-1}\mathcal R^{1/2})
\ge0$, which proves \eqref{extp:reciprocal-derivative}.
Since $P=0$ and $U=U_0$ on sufficiently deep levels, integration
gives \eqref{extp:integrated-U} and integrability for every $t<0$.

At $t=0$, Hopf's lemma gives $|\nabla\Phi|>0$ on $\partial D$,
so the level flow and the kernel variations extend smoothly from below.
The first variation yields
\[
c^*H(0)c=\int_{\partial D}
\frac{|\sum_jc_jk_{j,0}|^2}{|\nabla\Phi|}\,d\sigma>0
\qquad(c\ne0),
\]
where strictness follows from the maximum principle and independence
of the representers. Thus
$U(t)\to U(0)=\operatorname{tr}(M(0)H(0)^{-1})<\infty$.
Monotone convergence now gives
$B_\ell(0)=\lim_{t\uparrow0}B_\ell(t)\le U(0)-U_0<\infty$,
proving the endpoint assertion.
\end{proof}

\subsection{Proof of the stability estimate}

\begin{theorem}[Period stability]\label{extp:jet-stability}
If $t_0\le0$, $\varepsilon\ge0$, and $F_\ell(t_0)\le\varepsilon$,
then
\begin{equation}\label{extp:matrix-stability}
\int_{-\infty}^{t_0-\delta}\rho_\ell(s)\,ds
\le\frac{\varepsilon}{4\delta-2\varepsilon},
\qquad\delta>\varepsilon/2.
\end{equation}
The constant is independent of $\ell$.
\end{theorem}
\begin{proof}
We first prove the growth estimate
\begin{equation}\label{extp:derivative-obstruction}
F_\ell'(t)\ge\frac{4B_\ell(t)}{1+2B_\ell(t)},
\qquad -\infty<t<0.
\end{equation}
Put $v=F_\ell'(t)$ and $d_j=2j-\lambda_j(t)\ge0$.
Lemma~\ref{extp:speeds} gives $\sum_jd_j=v$.
If $v<2$, then
\begin{align*}
U(t)-U_0
&=\sum_{j=1}^\ell\frac{d_j}{2j(2j-d_j)}\\
&\le\frac{1}{2(2-v)}\sum_{j=1}^\ell d_j
=\frac{v}{2(2-v)}.
\end{align*}
Indeed, $d_j\le v$ and $2j\ge2$ bound each denominator below
by $2(2-v)$. Combining this with
\eqref{extp:integrated-U} and rearranging proves
\eqref{extp:derivative-obstruction}. If $v\ge2$, the same
inequality holds because $4B/(1+2B)<2$ for finite $B\ge0$.
The argument applies at critical parameters as well, since
$M$, $H$, and $F_\ell'$ are continuous there.

Let $B_*=B_\ell(t_0-\delta)$. The function $B_\ell$ is
nondecreasing, so for $t_0-\delta\le s<t_0$,
\[
F_\ell'(s)\ge\frac{4B_*}{1+2B_*}.
\]
Since $F_\ell(t_0-\delta)\ge0$, integration gives
$\varepsilon\ge4\delta B_* /(1+2B_*)$.
The assumption $\delta>\varepsilon/2$ permits rearrangement
to \eqref{extp:matrix-stability}. For $t_0=0$, pass to the
one-sided limit, using Proposition~\ref{extp:reciprocal}.
\end{proof}

This proves Theorem~\ref{main:stability}. Integration of
\eqref{extp:derivative-obstruction} also gives
\[
F_\ell(t)\ge\int_{-\infty}^t
\frac{4B_\ell(s)}{1+2B_\ell(s)}\,ds.
\]
For $\ell=1$, the coefficient is
$\rho_1=(K_t-K_t^e)/(K_t+K_t^e)$. For $\ell>1$, the
estimate controls the total period energy of an orthonormal basis
of the span of the jet representers.

\begin{corollary}[Total period energy]\label{period:total}
Under the hypotheses of Theorem~\ref{extp:jet-stability},
\[
\int_{-\infty}^{t_0-\delta}\operatorname{tr}P_t\,dt
\le\frac{\varepsilon}{2\delta-\varepsilon}.
\]
\end{corollary}
\begin{proof}
For $0\le x\le1$, $x\le2x/(2-x)$, so
$\operatorname{tr}P_t\le2\rho_\ell(t)$.
Apply \eqref{extp:matrix-stability}.
\end{proof}
At each regular level, let $e_1,\ldots,e_\ell$ be an orthonormal
basis of the span of the jet representers. Then
$\operatorname{tr}P_t=\sum_j\|(I-\Pi_e)e_j\|^2$.
Thus the integral is the accumulated total period energy over
$(-\infty,t_0-\delta]$, independently of the orthonormal basis
chosen at each regular level.

\section{Topology of the Green sublevels}\label{tp:section}
We combine the classical Green index formula
\cite[Section~2]{Saitoh} with the defect representation to locate
the threshold for positivity.

\subsection{The index formula}
Use the preceding notation for $D$, $q$, $\Phi$, and $D_t$.
At a critical point $p\ne q$ put
$\nu(p)=\operatorname{ord}_p(\partial_z\Phi)\ge1$.
Write $b_1$ for the first Betti number, which counts the holes
of a connected plane domain.

\begin{lemma}[Classical Green index formula]\label{tp:index}
There are finitely many critical points in $D\setminus\{q\}$.
For every regular $t<0$,
\begin{equation}\label{tp:index-formula}
b_1(D_t)=\sum_{\Phi(p)<t}\nu(p),
\qquad
b_1(D)=\sum_{\partial_z\Phi(p)=0}\nu(p).
\end{equation}
In particular, the number of holes cannot decrease along regular
Green sublevels.
\end{lemma}
\begin{proof}
Near $q$, the logarithmic term in $\Phi$ excludes critical points.
The Hopf lemma excludes them in a collar of $\partial D$.
The remaining zeros of the holomorphic function $\partial_z\Phi$
lie in a compact set and are isolated, hence finite.

Fix a regular $t$. Choose a smooth scalar function $\lambda$
that equals $|z-q|^2$ near $q$, is positive away from $q$, and
equals $1$ outside a small pole disk. Set $X=\lambda\nabla\Phi$
and identify planar vectors with complex numbers. Near $q$,
$\Phi=\log|z-q|+h$ with $h$ smooth and harmonic, so
\[
X(z)=(z-q)+|z-q|^2\nabla h(z),\qquad X(q)=0,\quad DX(q)=I_2.
\]
Thus $X$ extends smoothly across $q$ and
$\operatorname{ind}_qX=\operatorname{sgn}\det DX(q)=+1$.
At any other critical point $p$, write $\Phi=\operatorname{Re}W$
locally, so $W'(z)=2\partial_z\Phi(z)=(z-p)^{\nu(p)}a(z)$ with
$a(p)\ne0$. On a small counterclockwise circle about $p$, $W'$
has winding number $\nu(p)$. Since $\Phi_x+i\Phi_y=\overline{W'}$,
conjugation reverses this winding number, while $\lambda>0$
preserves it. Hence $\operatorname{ind}_pX=-\nu(p)$.

On $\partial D_t$, the outward unit normal is
$n=\nabla\Phi/|\nabla\Phi|$, so $X\cdot n=\lambda|\nabla\Phi|>0$.
Poincar\'e--Hopf on $\overline{D_t}$ and connectedness from
Lemma~\ref{lem:levels} therefore give
\[
1-b_1(D_t)=\chi(D_t)=1-\sum_{\Phi(p)<t}\nu(p).
\]
For $t$ sufficiently close to $0$, the boundary collar identifies
$D_t$ with $D$ up to diffeomorphism, proving the second formula.
\end{proof}

\subsection{The positivity threshold}

\begin{corollary}[Threshold for positivity]\label{tp:onset}\label{extp:jet-onset}
If $D$ is simply connected, then $F_\ell=0$ for every $\ell\ge1$.
Otherwise, let
\[
t_* =\min\{\Phi(p):p\in D\setminus\{q\},\ \partial_z\Phi(p)=0\}.
\]
For every $\ell\ge1$,
\[
F_\ell(t)=F_\ell'(t)=0\quad(t\le t_*),\qquad
F_\ell(t)>0,\quad F_\ell'(t)>0\quad(t_*<t<0).
\]
At every regular $t>t_*$, one has $E_\ell(t)>0$.
\end{corollary}
\begin{proof}
If $D$ is simply connected, the index formula gives no critical
points and every sublevel is simply connected. The conformal
formula in Lemma~\ref{mx:conformal-jets} gives $F_\ell=0$.
Otherwise, the same index formula shows that $t_*$ exists, that
sublevels below it are simply connected, and that regular sublevels
above it are multiply connected. Hence $F_\ell=F_\ell'=0$ below
$t_*$, and their continuity gives
$F_\ell(t_*)=F_\ell'(t_*)=0$.

Corollary~\ref{mx:rigidity} gives $E_\ell>0$ at every regular
level above $t_*$. By Theorem~\ref{mx:jet-defect},
\[
F_\ell'(t)=\int_{t_*}^tE_\ell(s)\,ds,
\qquad
F_\ell(t)=\int_{t_*}^t(t-s)E_\ell(s)\,ds.
\]
Every interval to the right of $t_*$ contains a regular subinterval,
where the density is positive. Both integrals are therefore
positive for $t>t_*$, including at later critical values.
\end{proof}

The threshold is independent of the jet order. For multiply
connected $D$, the curvature measure has relative support
$[t_*,0)$ in $(-\infty,0)$. The index formula counts subsequent
changes in the number of holes; positivity alone does not.

\section*{Disclosure of assistance}
An OpenAI language model assisted with proof development, literature
comparison, and editorial revision.

\end{document}